\documentclass{amsart}

\usepackage{mystyle}

\makeatletter
\def\@tocline#1#2#3#4#5#6#7{\relax
  \ifnum #1>\c@tocdepth 
  \else
    \par \addpenalty\@secpenalty\addvspace{#2}%
    \begingroup \hyphenpenalty\@M
    \@ifempty{#4}{%
      \@tempdima\csname r@tocindent\number#1\endcsname\relax
    }{%
      \@tempdima#4\relax
    }%
    \parindent\z@ \leftskip#3\relax \advance\leftskip\@tempdima\relax
    \rightskip\@pnumwidth plus4em \parfillskip-\@pnumwidth
    #5\leavevmode\hskip-\@tempdima
      \ifcase #1
       \or\or \hskip 1em \or \hskip 2em \else \hskip 3em \fi%
      #6\nobreak\relax
    \dotfill\hbox to\@pnumwidth{\@tocpagenum{#7}}\par
    \nobreak
    \endgroup
  \fi}
\makeatother

\begin{document}

\title{The Motivic cofiber of $\tau$}

\author{Bogdan Gheorghe}
\email{gheorghebg@wayne.edu}

\address{Department of Mathematics\\ Wayne State University\\
Detroit, MI 48202, USA}
\date{\today}


\begin{abstract}
Consider the Tate twist $\tau \in H^{0,1}(\Ss)$ in the mod 2 cohomology of the motivic sphere. After 2-completion, the motivic Adams spectral sequence realizes this element as a map $\tau \colon S^{0,-1} \lto \Ss$, with cofiber $\Ct$. We show that this motivic 2-cell complex can be endowed with a unique $E_{\infty}$ ring structure. Moreover, this promotes the known isomorphism $\piss \Ct \cong \Ext^{\ast,\ast}_{BP_{\ast}BP}(BP_{\ast},BP_{\ast})$ to an isomorphism of rings which also preserves higher products. 

We then consider the closed symmetric monoidal category $(\CtMod, - \smasCt -)$ which lives in the kernel of Betti realization. Given a motivic spectrum $X$, the $\Ct$-induced spectrum $X \smas \Ct$ is usually better behaved and easier to understand than $X$ itself. We specifically illustrate this concept in the examples of the mod 2 Eilenberg-Maclane spectrum $\HF$, the mod 2 Moore spectrum $\Ss/2$ and the connective hermitian $K$-theory spectrum $\kq$.
\end{abstract}

\maketitle


\tableofcontents

\section{Introduction}

\subsection{The Setting}

The mod 2 cohomology of the motivic sphere spectrum $\Ss$ over $\Spec \C$ was computed by Voevodsky in \cite{VoeZ2}, and is given by
\begin{equation*}
\HFco(\Ss) \cong \F_2[\tau] \qquad \text{ where } |\tau| = (0,1).
\end{equation*}
Denote the mod 2 motivic Steenrod algebra of operations $\left[ \HF, \HF  \right]_{\ast,\ast}$ by $\AC$. One can run the motivic Adams spectral sequence 
\begin{equation*}
\Ext_{\AC}( \F_2[\tau], \F_2[\tau]) \Longrightarrow \piss\left( (\Ss)^{\smas}_2 \right),
\end{equation*}
as constructed in \cite{MorelAdams}, \cite{DuggerIsaksenMASS}, \cite{HKOAdams}. Observe that the $E_2$-page contains a non-trivial element in Adams filtration 0, namely multiplication by $\tau$ on $\F_2[\tau]$. This is different from the topological Adams spectral sequence for $S^0$, where the only elements in Adams filtration 0 are the identity map and the zero map. It is easy to see that this element survives to the $E_{\infty}$-page as it cannot be involved with any differential for degree reasons. Therefore, it detects a map
\begin{equation*}
S^{0,-1} \stackrel{\tau}{\lto} (\Ss)^{\smas}_2,
\end{equation*}
whose Hurewicz image is $\tau \in \HFho( (\Ss)^{\smas}_2)$. 

To avoid complications about the existence of a non-completed version of this map, we will now work 2-completed. Recall that 2-completion is given by the $E$-Bousfield localization at either the Moore spectrum $\Ss/2$ or the Eilenberg-Maclane spectrum $\HF$. In particular, the 2-completed sphere $L_E \Ss$ is also an $E_{\infty}$ ring spectrum and admits a good category of ($2$-completed) modules. We will from now on work in the 2-completed category, i.e., in modules over the 2-completed sphere. We will denote the 2-completed sphere and the smash product in 2-completed spectra simply by $\Ss$ and $- \smas -$. With this notation, the motivic Adams spectral sequence produces a non-trivial map $S^{0,-1} \stackrel{\tau}{\lto} \Ss$. 

Recall that the Betti realization functor $\Betti$ goes from (here 2-completed) motivic spectra $\SptC$ over $\Spec \C$ to classical (2-completed) spectra $\Spt$. This functor is for example constructed in \cite[2.6]{DuggerIsaksenMASS}, \cite[Appendix A.7]{PPR} or \cite[Chapter 4]{Joa}, and is induced by taking $\C$-points of the involved $\C$-schemes. It is a left adjoint, with right adjoint usually denoted $\Sing$, and admits the constant functor $c$ as a section \cite{Levine14}. The situation is summarized in the diagram
\begin{center}
\begin{tikzpicture}
\matrix (m) [matrix of math nodes, row sep=3em, column sep=4em]
{ \SptC & \Spt. \\};
\path[thick, -stealth, font=\small]
(m-1-1.10) edge node[above] {$ \Betti $} (m-1-2.170)
(m-1-2.192) edge node[below] {$ \Sing $} (m-1-1.350);
\path[thick, dashed, -stealth, font=\small]
(m-1-2.135) edge[bend right=45] node[above] {$ c $} (m-1-1);
\end{tikzpicture}
\end{center}
The Betti realization functor $\Betti$ therefore induces a split-surjection
\begin{equation*}
\pisw(\Ss) \surj \pi_s(S^0),
\end{equation*}
with section induced by the constant functor $c$. Moreover, it sends the map $S^{0,-1} \stackrel{\tau}{\lto} \Ss$ to the identity $S^0 \stackrel{\id}{\lto} S^0$, as shown in \cite[Section 2.6]{DuggerIsaksenMASS}. Computationally, the Betti realization functor $\Betti$ can thus be interpreted as sending the element $\tau$ to $1$. For example, on the homotopy of the mod 2 Eilenberg-Maclane spectrum it induces the quotient map 
\begin{equation*}
\piss(\HF) \cong \M \surj \pi_{\ast}(\HF) \cong \F_2,
\end{equation*}
which imposes the relation $\tau = 1$. Observe that there is another surjection $\M \surj \F_2$ with same source and target, namely the quotient map imposing the relation $\tau = 0$. One can thus ask if this map is also induced by a functor between $\SptC$ and another homotopy theory. To answer this question, we are led to study the homotopy theoretical analogue of the algebraic operation of setting $\tau = 0$, which is to take the cofiber of the map $\tau$. Consider the cofiber sequence
\begin{equation} \label{intro:Ctdef}
S^{0,-1} \stackrel{\tau}{\lto} \Ss \lto \Ct \lto S^{1,-1},
\end{equation}
where we denote the cofiber of the map $\tau$ by $C\tau$. This 2-cell complex already appeared in \cite{StableStems}, where it is studied via its motivic Adams-Novikov spectral sequence. More precisely, it is proven that its Adams-Novikov spectral sequence collapses at the $E_2$-page with no possible hidden extensions. This provides a surprising isomorphism 
\begin{equation} \label{intro:CtBP}
\Ext^{\ast,\ast}_{BP_{\ast}BP}(BP_{\ast}, BP_{\ast}) \cong \piss(\Ct),
\end{equation}
connecting two objects which are a priori unrelated. The left hand side is the cohomology of the classical (non-motivic) Hopf algebroid $(BP_{\ast},BP_{\ast}BP)$ and is very important in chromatic homotopy theory. In particular, it is the $E_2$-page of the Adams-Novikov spectral sequence for the topological sphere $S^0$. Notice that since it is the cohomology of a dga, namely the cobar complex associated to $(BP_{\ast},BP_{\ast}BP)$, it admits products and higher Massey products. All this algebraic structure gets transferred to the motivic homotopy groups $\piss(\Ct)$, formally endowing it with a (higher) ring structure. One can thus ask if this algebraic ring structure can be lifted to a topological ring structure on $\Ct$. The first goal of this paper is to answer this question, which we do in Section \ref{sec:CtEinfty} by the following results.

\begin{thm} \label{intro:thmCtEoo}
There exists a unique $E_{\infty}$ ring structure on $\Ct$. 
\end{thm}

We now explain what we mean by a motivic $E_{\infty}$ ring spectrum, and refer to Section \ref{sec:motoperads} for more details. Since $\SptC$ is enriched over simplicial sets, one can talk about algebras over operads in simplicial sets. In fact, operads in simplical sets embedded in the motivic world are sometimes called constant operads. In this paper, we say that a motivic spectrum admits an $E_{\infty}$ ring structure if it admits an algebra structure over a constant $E_{\infty}$ operad, i.e., over any usual $E_{\infty}$ operad in simplicial sets. We warn the reader that similarly to the equivariant case of \cite{BlumbergHill}, this notion of motivic $E_{\infty}$ ring spectra is probably not the same as strictly commutative algebras in $\SptC$.

There are two main tools involved in proving Theorem \ref{intro:thmCtEoo}. We first use elementary techniques with triangulated categories to produce a unital, associative and commutative monoid in the homotopy category $\Ho(\SptC)$. We then rigidify this ring structure using Robinson's $E_{\infty}$ obstruction theory \cite{RobinsonGamma}. By tracing back to the origin of the isomorphism \eqref{intro:CtBP}, we can now show that the algebraic structure on $\piss(\Ct)$ does come from $\Ct$.

\begin{prop} \label{intro:CtEoo}
The isomorphism \eqref{intro:CtBP}
\begin{equation*}
\piss(\Ct) \cong \Ext^{\ast,\ast}_{BP_{\ast}BP}(BP_{\ast}, BP_{\ast})
\end{equation*}
is an isomorphism of rings which sends Toda brackets in $\piss$ to Massey products in $\Ext$, and vice-versa.
\end{prop}

Let's point out that the additive version of this theorem was already exploited by Isaksen in \cite{StableStems} to gain knowledge about the classical Adams-Novikov $E_2$-page. The idea is to compute $\piss(\Ct)$ in a range using its motivic Adams spectral sequence and the knowledge of $\piss(\Ss)$ in this range. Having a multiplicative structure available improves the correspondence in an obvious manner.

Having considered the cofiber of multiplication by $\tau$, one can look at the less severe quotients $\quotient{\Ss}{\tau^n} \eqqcolon C\tau^n$ of multiplication by $\tau^n$. Using similar techniques as in Theorem \ref{intro:CtEoo}, one can show that every spectrum $C\tau^n$ is uniquely an $A_{\infty}$ ring spectrum, and that it is homotopy commutative. Our method does not apply to show that $\Ctn$ is $E_{\infty}$ for any $n$, as the necessary obstruction groups do not vanish. It is however probably the case that every $\Ctn$ is also an $E_{\infty}$ ring spectrum, since $C\tau^n$ gets closer and closer to $\Ss$ as $n$ increases. One can also show that for any $n$, the natural reduction map $\Ctn \lto \Ctnn$ is a ring map. All together, these spectra sit in a tower 
\begin{equation} \label{intro:diagtowerofCtn}
\begin{tikzpicture}
\matrix (m) [matrix of math nodes, row sep=2em, column sep=3em]
{  \Ct & C\tau^2           & C\tau^3            & C\tau^4 & \cdots & (\Ss)^{\smas}_{\tau}  \\
       & \Sigma^{0,-1} \Ct & \Sigma^{0,-2} \Ct  & \Sigma^{0,-3} \Ct & & \\};
\path[thick, -stealth]
(m-1-2) edge (m-1-1)
(m-1-3) edge (m-1-2)
(m-1-4) edge (m-1-3)
(m-1-5.171) edge (m-1-4)
(m-1-6) edge (m-1-5.5)

(m-2-2) edge  (m-1-2)
(m-2-3) edge (m-1-3)
(m-2-4) edge (m-1-4);
\end{tikzpicture}
\end{equation} 
of ring spectra and ring spectra maps, where each layer is a copy of $\Ct$. Since the completion map $\Ss \lto (\Ss)^{\smas}_{\tau}$ induces an isomorphism on homotopy groups, this tower reconstructs the sphere spectrum $\Ss$ by increasing the exponent of $\tau$-torsion. More precisely, one can consider the $\tau$-Bockstein spectral sequence for $\Ss$, which is the spectral sequence induced by applying $\piss$ to the tower \eqref{intro:diagtowerofCtn}. Surprisingly, this spectral sequence contains the same information as the motivic Adams-Novikov spectral sequence computing $\piss(\Ss)$. In fact, the $E_1$-page of the $\tau$-Bockstein spectral sequence is isomorphic to the $E_2$-page of the motivic Adams-Novikov spectral sequence. Moreover, by \cite[Lemma 15]{HKO11}, the motivic Adams-Novikov spectral sequence has only odd differentials, which are all of the form $d_{2r+1}(x) = \tau^r y$. Such a differential corresponds to a $d_r$ differential of the $\tau$-Bockstein spectral sequence, giving a one-to-one correspondence between the differentials of each spectral sequence. This implies that the $E_{2r+2}$ page of the motivic Adams-Novikov spectral sequence is isomorphic to the $E_{r+1}$ page of the $\tau$-Bockstein spectral sequence. 

\vspace{0.3cm}

With an $E_{\infty}$ ring structure in hand, any good model for motivic spectra produces a closed symmetric monoidal category of $\Ct$-modules with the relative smash product $- \smasCt -$, and a free-forget adjunction 
\begin{equation} \label{eq:adjSModCtMod}
\SptC \stackrel{-\smas \Ct}{\adj} \CtMod.
\end{equation}
The remainder of this paper is devoted to the task of better understanding the category $\CtMod$. In Lemma \ref{lem:BettiRelofCtMod} we show that the Betti realization of any $\Ct$-module is contractible, which means that the category of $\Ct$-modules lies in the kernel of Betti realization. This does not mean that the motivic spectrum $\Ct$ does not have topological applications, as there are other bridges between motivic and classical homotopy theory. Such a bridge is for example given by Proposition \ref{intro:CtEoo}, relating the homotopy groups of the motivic spectrum $\Ct$ with the cohomology of the Hopf algebroid $(BP_{\ast},BP_{\ast}BP)$. 

One strength of the category $\CtMod$ is that it is relatively easy to work with $\Ct$-modules. One first observes this phenomenon during the process of proving that $\Ct$ admits an $E_{\infty}$ ring structure, with the many obstruction groups vanishing for degree reasons. We observe a similar phenomenon with related motivic spectra. For example, we can completely describe the ring spectra $\Ct \smas \Ct$ and $\End(\Ct)$ by using elementary techniques. In joint work with Zhouli Xu and Guozhen Wang \cite{GWX}, we provide an equivalence between some category of (cellular) $\Ct$-modules and some category of derived $BP_{\ast}BP$-comodules. In particular, this implies that the homotopy category of cellular $\Ct$-modules is algebraic in the sense of \cite{Schwede}. This is another reason why it feels easier to manipulate motivic spectra living in $\CtMod$, since algebraic categories are usually better behaved than topological categories. For example, algebraic categories admit a $\mathcal{D}(\Z)$-enrichment which implies many pleasant properties. We refer the reader to Remark \ref{rem:exwhyalgcatarebetter} for a concrete such example, and to \cite{Schwede} for more details. 

The category  $\CtMod$ is the universal place in which the element $\tau$ has been killed. The strength of this benign statement lies in the fact that many motivic spectra naturally land in $\CtMod$, since at some point we were led to mod out by $\tau$ for one reason or another. For example, the relation $\tau \eta^4 = 0 \in \piss$ implies that the $\eta$-inverted sphere $\Ss[\eta^{- 1}]$, computed in \cite{GuiI} and \cite{AM14}, lives in $\CtMod$. More generally, one can show that any element $x \in \pisw$ with $s \neq 0$ admits a relation of the type $\tau^a x^b = 0$. 
In particular, inverting any such non-nilpotent element $x$ when $a=1$ yields a spectrum that naturally lives in $\CtMod$. In particular, this phenomenon applies to the exotic Morava $K$-theories $K(w_n)$ detecting the motivic $w_n$-periodicity that will appear in a future paper \cite{GheKwn}.


\vspace{0.3cm}

We now describe the last Section of this paper, where we explicitly compute the homotopy of some specific $\Ct$-modules induced through the adjunction \eqref{eq:adjSModCtMod}. Given a spectrum $X$, we call the induced $\Ct$-module $X \smas \Ct$ a $\Ct$-induced spectrum.

One of the first spectra to understand in $\CtMod$ is the $\Ct$-induced mod 2 Eilenberg-Maclane spectrum. This is the spectrum $\HF \smas \Ct$ that we denote by $\HCt$, which we will treat as a cohomology theory. Given any $\Ct$-module $X$ one can consider the $\Ct$-linear mod 2 (co)homology of $X$ defined by the homotopy of the spectra
\begin{equation*}
\FCt(X, \HCt) \qquad \text{ and } \qquad X \smasCt \HCt.
\end{equation*}
Here $\FCt(-,-)$ denotes the $\Ct$-linear function spectrum and $- \smasCt -$ denotes the relative smash product in $\CtMod$. Observe for example that since
\begin{equation*}
X \smasCt \HCt = X \smasCt (\HF \smas \Ct) \simeq X \smasCt \Ct \smas \HF \simeq X \smas \HF,
\end{equation*}
the $\Ct$-linear $\HCt$-homology of $X$ is isomorphic to the $\HF$-homology of the underlying spectrum of $X$. Consider the $\Ct$-linear (Steenrod) algebra of $\HCt$-(co-)operations, given by the homotopy of the spectra
\begin{equation*}
\FCt(\HCt, \HCt) \qquad \text{ and } \qquad \HCt \smasCt \HCt.
\end{equation*}
These are the relevant (co-)operations acting on the $\Ct$-linear (co)homology of $\Ct$-modules, which we compute in Section \ref{sec:HFtwomodtau}. Recall that the dual mod 2 motivic Steenrod algebra over $\Spec \C$ is given by
\begin{equation*}
\piss( \HF \smas \HF) \cong \quotient{ \F_2[\tau][\xi_1, \xi_2, \ldots,\tau_0, \tau_1, \ldots]}{\tau_i^2 = \tau \xi_{i+1}}.
\end{equation*}
See Section \ref{sec:SteenrodAlg} for more details. The following computation follows easily.

\begin{prop} \label{intro:HF2CtSteenrod}
The Hopf algebra of $\Ct$-linear co-operations of $\HF$ is given by
$$\F_2[\xi_1, \xi_2, \ldots] \otimes E(\tau_0, \tau_1, \ldots).$$
\end{prop}

One can also as usual consider $\HCt$ as a (co)homology theory on $\SptC$, and define the $\HCt$-homology and cohomology of any motivic spectrum $X$ by the homotopy of the spectra
\begin{equation*}
F(X, \HCt) \qquad \text{ and } \qquad X \smas \HCt.
\end{equation*}
The associated (co-)operations acting on the $\HCt$-(co-)homology of any spectrum is given by the homotopy of the spectra
\begin{equation*}
F(\HCt, \HCt) \qquad \text{ and } \qquad \HCt \smas \HCt.
\end{equation*}
We also compute these Hopf algebras in Section \ref{sec:HFtwomodtau}, in particular we get the following.

\begin{prop}
The Hopf algebra of co-operations of $\HF$ is given by 
$$\F_2[\xi_1, \xi_2, \ldots] \otimes E(\tau_0, \tau_1, \ldots) \otimes E(\btau).$$
\end{prop}

The extra co-operation $\btau$ is primitive in the coalgebra structure. It is induced by the $\tau$-Bockstein corresponding to the composite 
\begin{equation*}
\Ct \stackrel{p}{\lto} S^{1,-1} \stackrel{i}{\lto} \SCt
\end{equation*}
of the projection of $\Ct$ on its top cell $S^{1,-1}$, followed by the inclusion as its bottom cell. It also appears in the Hopf algebra of operations, where it really deserves its name of $\tau$-Bockstein. More precisely, the operation $\btau$ in cohomology is to $\tau$ as the usual Bockstein $\beta = \Sq$ is to the element 2. In particular, it allows to reconstruct $\HF$ via a motivic analogue of the Postnikov tower that runs in the weight direction. At each stage, the layer is a copy of $\HCt$, and the boundary map composed with the $k$-invariant is given by $\btau$. However, if one is only interested in studying $\Ct$-modules internally to the category $\CtMod$, this is a noisy element and one should use $\Ct$-linear (co)homology.


\vspace{0.3cm}

Given a spectrum $X \in \SptC$, the homotopy groups of the $\Ct$-induced spectrum $X \smas \Ct$ are an extension of the $\tau$-torsion and the residue mod $\tau$ of $\piss(X)$. Even though this proves to be very complicated in general, one principle that appears is the following. If all the obstructions to the spectrum $X \in \SptC$ possessing some property or structure are $\tau$-torsion, then the $\Ct$-induced spectrum $X \smas \Ct$ posses the desired property or structure. Here are a few such examples that we study in Section \ref{sec:Ctmodules}.

\vspace{0.1cm}

Start with the 2-completed motivic mod 2 Moore spectrum $\Ss/2$. The Toda bracket $\langle 2,\eta,2 \rangle \ni \tau \eta^2$ is the obstruction to both endowing it with a left unital multiplication, and to a $v_1^1$-self map. We can prove the following results about the $\Ct$-induced Moore spectrum, which we denote by $\Stwotau$.

\begin{thm}
The $\Ct$-induced motivic mod 2 Moore spectrum $\Stwotau$ admits a unique structure of an $E_{\infty}$ $\Ct$-algebra.
\end{thm}

\begin{prop}
The $\Ct$-induced motivic mod 2 Moore spectrum $\Stwotau$ admits a $v_1^1$-self map
\begin{equation*}
\Sigma^{2,1} \Stwotau \stackrel{v_1}{\lto} \Stwotau.
\end{equation*}
\end{prop}

\vspace{0.1cm}

We also study the 2-completed connective algebraic and hermitian $K$-theory spectra $\kgl$ and $\kq$. Denote their $\Ct$-induced spectra by $\kglCt \coloneqq \kgl \smas \Ct$ and $\kqCt \coloneqq \kq \smas \Ct$. The case of algebraic $K$-theory is very simple as its invariant are $\tau$-free. For the case of the connective hermitian $K$-theory spectrum (constructed over $\Spec \C$ in \cite{IS}, where it is denoted $ko$) we prove the following result.

\begin{prop}
The $\Ct$-induced connective hermitian $K$-theory spectrum $\kqCt$ has homotopy groups
\begin{equation*}
\piss( \kqCt) \cong \quotient{\hat{\Z}_2[v_1^2,\eta]}{2 \eta}.
\end{equation*}
\end{prop}

Recall that the homotopy of the motivic spectrum $\kq$ contains the 8-fold Bott periodicity element $v_1^4$, but does not contain $v_1^2$. One can compute the homotopy of $\kq$ and $\kqCt$ from their cohomology via the motivic May spectral sequence, followed by the motivic Adams spectral sequence. In the case of $\kq$, there is a motivic May differential supported by $v_1^2$ and with $\tau$-torsion target. In the case of $\kqCt$, the $\tau$-torsion target of this differential gets shifted. We then resolve a hidden extension of this shifted element to show that it is in fact the periodicity element $v_1^2$. More precisely, we show that this element is a square root of the usual 8-fold Bott periodicity, making the $\Ct$-induced spectrum $\kqCt$ 4-fold periodic. In chromatic motivic language, up to the $v_0$-extensions this can be rewritten as $\quotient{\F_2[v_0, v_1^2, w_0]}{v_0 w_0}$. The relation $v_0 w_0 = 0$ is clear as it is already present in $\piss(\Ss)$, but this shows that $v_1^2$ and $w_0$ can coexist without any relation between them.

\subsection{The Choice of Prime $p = 2$} 

This paper is written in a $p$-completed setting, where we chose the prime $p=2$. However, the main results also apply to odd primes. In short, the $H\F_p$-based motivic Adams spectral sequence produces the map
\begin{equation*}
S^{0,-1} \stackrel{\tau}{\lto} (\Ss)^{\smas}_p,
\end{equation*}
after $p$-completing the target for any prime $p$. Denote again its cofiber by $\Ct$, where the prime $p$ does not appear in the notation. The isomorphism
\begin{equation*}
\Ext^{\ast,\ast}_{BP_{\ast}BP}(BP_{\ast}, BP_{\ast}) \cong \piss(\Ct)
\end{equation*}
still holds for any prime, producing the same vanishing regions in the homotopy of $\Ct$, and thus endowing $\Ct$ with an $E_{\infty}$ ring structure. 

For odd primes $p$, the motivic story is somehow easier since it is more closely related to the classical story. In particular, in the case of odd primes, the motivic Steenrod algebra (and its dual) are isomorphic as Hopf algebras to the classical Steenrod algebra (and its dual), adjoint a formal variable $\tau$. This is not the case for $p=2$, for example because of the relation $\tau_i^2 = \tau \xi_{i+1}$ in the dual motivic Steenrod algebra.

\subsection{Organization} \mbox{}

\noindent \textbf{Section \ref{sec:BackgroundandNot}.} This Section contains a brief summary of the motivic homotopy theory needed in order to define the spectrum $\Ct$. This contains a recall of the motivic category of spectra over $\Spec \C$, some functors relating $\SptC$ with the topological category $\Spt$, the mod 2 motivic cohomology, the structure of the mod 2 motivic Steenrod algebra and its dual, and the motivic Adams spectral sequence. After introducing the spectrum $\Ct$, we explain some vanishing regions both in its homotopy groups $\piss(\Ct)$ and in the homotopy classes of self-maps $\left[ \Ct, \Ct \right]_{\ast,\ast}$. These results will be mostly used to endow $\Ct$ with an $E_{\infty}$ ring structure.

\vspace{0.3cm}

\noindent \textbf{Section \ref{sec:CtEinfty}.} We first explain the notion of motivic $A_{\infty}$ and $E_{\infty}$ ring spectra that we will use in this paper, and adapt Robinson's obstruction theory \cite{RobinsonGamma} to the motivic setting. We then apply this obstruction theory to endow the spectrum $\Ct$ with an $E_{\infty}$ ring structure.

\vspace{0.3cm}

\noindent \textbf{Section \ref{sec:computationsCt}.} In this Section we compute the homotopy types of the $E_{\infty}$ ring spectrum $\Ct \smas \Ct$ and of the $A_{\infty}$ ring spectrum $\ECt$.

\vspace{0.3cm}

\noindent \textbf{Section \ref{sec:Ctmodules}.} This Section is about the symmetric monoidal category $\CtMod$. We start by showing some generalities on $\Ct$-modules. We then analyze more precisely a few specific $\Ct$-induced spectra:
\begin{enumerate}
\item We compute the Steenrod algebra of operations and co-operations on the $\Ct$-induced mod 2 Eilenberg-Maclane spectrum $\HF \smas \Ct$.
\item We show that the $\Ct$-induced mod 2 Moore spectrum $\Stwotau$ admits a unique $E_{\infty}$ structure as a $\Ct$-algebra, and that it admits a $v_1^1$-self map.
\item We compute the homotopy groups of the $\Ct$-induced connective algebraic and hermitian $K$-theories $\kgl \smas \Ct$ and  $\kq \smas \Ct$. In particular, a hidden extension shows that $\kq \smas \Ct$ contains a 4-fold periodicity by the element $v_1^2$, which is the square root of the usual 8-fold Bott periodicity observed in $\kq$.
\end{enumerate}

\subsection{Acknowledgment} 

The author is grateful for contributions from Dan Isaksen, Nicolas Ricka, Prasit Bhattacharya, Sean Tilson, Peter May and Mike Hill.

\section{Notation and Background on $\Ct$} \label{sec:BackgroundandNot}

In this Section we give some brief background on motivic homotopy theory over $\Spec \C$ as well as properly introduce the spectrum $\Ct$. For a more detailed introduction to motivic homotopy theory we refer the reader to \cite{MorelA1}, \cite{MorelVoevodsky}. Most of our notation agrees and is taken from \cite{StableStems}.


\subsection{Motivic Spaces and Spectra over $\Spec \C$}

Denote by $\SpcC$ the category of \emph{(pointed) motivic spaces} over $\Spec \C$ as defined in \cite{MorelVoevodsky}. This category admits a well-behaved homotopy theory, as it supports a closed symmetric monoidal, proper, simplicial and cellular model structure. The paper \cite[Chapter 2]{Pelaez} is a good source for a careful construction of these model structures. There is a realization functor 
\begin{equation*}
\SpcC \stackrel{\Betti}{\ltooo} \Top,
\end{equation*}
from motivic spaces over $\Spec \C$ to topological spaces called \emph{Betti realization}. This functor is for example constructed in \cite[2.6]{DuggerIsaksenMASS}, \cite[Appendix A.7]{PPR} or \cite[Chapter 4]{Joa}, and is induced by taking $\C$-points of the involved $\C$-schemes. It is a strict symmetric monoidal left Quillen functor, whose right adjoint is usually denoted by $\Sing$. In the same spirit as equivariant homotopy theory, motivic homotopy theory has two different types of spheres. We will denote the 1-dimensional simplicial sphere by $S^{1,0} \in \SptC$ and the geometric sphere $\G_m$ by $S^{1,1} \in \SptC$. The first coordinate $m$ in the notation $S^{m,n}$ indicates the \emph{topological dimension} of the sphere, and it is not hard to see that it Betti realizes to the topological sphere $S^m$. The second coordinate $n$ indicates the \emph{weight}, or the \emph{Tate twist} of the sphere $S^{m,n}$. Over $\Spec \C$, the projective line $\PP^1$ is a 2-dimensional topological sphere, whose homotopy type is described by the equation
\begin{equation} \label{eq:P1}
\PP^1 \simeq S^{1,0} \smas S^{1,1} \simeq S^{2,1}.
\end{equation}
The category of \emph{motivic $(\PP^1$-$)$spectra} $\SptC$ over $\Spec \C$ is constructed by stabilizing with respect to the sphere $\PP^1$, i.e., inverting the functor $- \smas \PP^1$. Observe that equation \eqref{eq:P1} implies that this is equivalent to inverting smashing with both fundamental spheres $- \smas S^{1,0}$ and $- \smas S^{1,1}$. This provides a bigraded suspension functor that we denote by $\Sigma^{m,n} = - \smas S^{m,n}$. Smashing with the simplicial sphere $\Sigma = \Sigma^{1,0} = - \smas S^{1,0}$ corresponds to the shift functor of the triangulated structure on the homotopy category. The category of motivic spectra $\SptC$ also supports good model structures which are closed symmetric monoidal with respect to the smash product $- \smas -$, proper, simplicial and cellular. The paper \cite[Chapter 2]{Pelaez} constructs these models in details. Moreover, the realization and singular pair stabilizes to a Quillen adjunction\footnote{Since the Betti realization of $\PP^1$ is the topological sphere $\PP^1(\C) \simeq S^2$, taking $\C$-points lands in the category of $S^2$-spectra, i.e., spectra with bonding maps $S^2 \smas X_n \lto X_{n+1}$. This is also a model for stable homotopy theory, see \cite[Section 4.1]{Joa} for more details.}
\begin{equation*}
\SptC \overunder{\Betti}{\Sing}{\adj} \Spt,
\end{equation*}
where the Betti functor $\Betti$ is strict symmetric monoidal, see for example \cite[A.45]{PPR}.

Given two spectra $X,Y \in \SptC$, the closed symmetric monoidal structure provides a \emph{function motivic spectrum} that we denote by $F(X,Y) \in \SptC$. When $X=Y$, we will usually write $\End(X) = F(X,X)$. As usual, we will denote the abelian group of homotopy classes of maps between $X$ and $Y$ by $[X,Y]$. When the source spectrum is a sphere $X = S^{s,w}$, the abelian group
\begin{equation*}
\pi_{s,w} (Y) \coloneqq [S^{s,w}, Y]
\end{equation*}
is called the homotopy group of $Y$ in \emph{stem $s$} and \emph{weight $w$}. The relation between the two is given by the usual adjunction between the smash product and the function spectrum. After taking homotopy, this becomes the equation
\begin{equation*}
\pi_{s,w} (F(X,Y)) \cong [ \Sigma^{s,w} X, Y].
\end{equation*}


\subsection{The Motivic Steenrod Algebra and the Adams Spectral Sequence} \label{sec:SteenrodAlg}

Denote by $H\Z$ Voevodsky's motivic Eilenberg-Maclane spectrum representing integral motivic cohomology on schemes \cite[Section 6.1]{Voe}. Denote by $\HF$ the cofiber of multiplication by 2 on $H\Z$, which sits in the cofiber sequence
\begin{equation*}
H\Z \stackrel{\cdot 2}{\ltoo} H\Z \lto \HF.
\end{equation*}
The spectrum $\HF$ represents mod 2 motivic cohomology on schemes. The coefficients of this spectrum were computed in \cite{VoeZ2} and are given by
\begin{equation*}
\HFho (\Ss) \cong \Ftwo[\tau] \qquad \text{ for } |\tau| = (0,-1).
\end{equation*}
Dually, the motivic cohomology of a point is 
\begin{equation*}
\HFco(\Ss) \cong \Ftwo[\tau] \qquad \text{ for } |\tau| = (0,1),
\end{equation*}
where we abuse notation and use the same symbol $\tau$ to denote the Tate twist element in homology and its dual in cohomology. We use the same notation as in \cite{StableStems} for the coefficients
\begin{equation*}
\M \coloneqq \HFco (\Ss) \cong \Ftwo[\tau] \qquad \text{ and } \qquad \Mdual \coloneqq \HFho (\Ss) \cong \Ftwo[\tau].
\end{equation*}
We write $\AC$ for the mod 2 motivic Steenrod algebra, i.e., the ring of stable cohomology operations on the motivic spectrum $\HF$. Its structure has been computed by Voevodsky in \cite{Voered}, \cite{VoeEM} : it is the bigraded Hopf algebra over $\M$ given by 
\begin{equation*}
\AC \cong \quotient{\M \langle \Sq, \Sqq, \ldots \rangle}{\text{Adem relations}}.
\end{equation*}
Observe that as in topology, it is generated by the Steenrod squares $\Sqnone$ with the Adem relations between them. The Tate twist $\tau \in \M$ has bidegree $|\tau|=(0,1)$, and the Steenrod squares have bidegrees $|\Sqn| = (2n,n)$ and $|\Sqnn| = (2n+1,n)$. Since we work at $p=2$, the first square $\Sq = \beta$ is again the usual Bockstein operation coming from the short exact sequence of abelian groups
\begin{equation*}
0 \lto \Z/2 \stackrel{\cdot 2}{\lto} \Z/4 \lto \Z/2 \lto 0.
\end{equation*}
The dual motivic Steenrod algebra 
\begin{equation} \label{eq:dualSteenrodalg}
\ACdual \cong  \quotient{\Mdual[\xi_1, \xi_2, \ldots, \tau_0, \tau_1, \ldots]}{\tau \xi_{i+1} = \tau_{i}^2},
\end{equation}
was also computed by Voevodsky in \cite{Voered}. Because we are now in homology, the Tate twist $\tau \in \Mdual$ has bidegree $|\tau|=(0,-1)$. The $\xi_i$'s and $\tau_i$'s have bidegrees $|\xi_i| = (2^{i+1} - 2,2^i -1)$ and $|\tau_i| =(2^{i+1}-1,2^i-1)$. The coproduct is given by the formulas
\begin{equation*}
\Delta(\xi_n) = \sum \xi_{n-k}^{2^k} \otimes \xi_k \qquad \text{ and } \qquad \Delta(\tau_n) = \tau_n \otimes 1 + \sum \xi_{n-k}^{2^k} \otimes \tau_k
\end{equation*}
One can now run the motivic Adams spectral sequence 
\begin{equation*}
\Ext_{\AC}(\M,\M) \Longrightarrow \piss ((\Ss)^{\smas}_2)
\end{equation*}
constructed in \cite{MorelAdams}, \cite{DuggerIsaksenMASS}, \cite{HKOAdams}, that converges to the homotopy groups of the 2-completed motivic sphere $(\Ss)^{\smas}_2$. The $\AC$-module map
\begin{equation*}
\M \stackrel{\cdot \tau}{\lto} \M
\end{equation*}
is an element in $\Hom = \Ext^0$ of Adams filtration 0 as $\tau$ is central in $\AC$. This element survives to the $E_{\infty}$-page as it cannot be involved with any differential for degree reasons. Therefore, it detects a map
\begin{equation*}
S^{0,-1} \stackrel{\tau}{\lto} (\Ss)^{\smas}_2,
\end{equation*}
whose Hurewicz image is the element $\tau \in \HFho ( (\Ss)^{\smas}_2)$.

It is crucial for us that this element $\tau$ exists in the homotopy groups of the motivic sphere spectrum, and thus acts on the homotopy of any motivic spectrum. To avoid complications with the lift of this element to the non-completed sphere, we will now work 2-completed. Recall that 2-completion is given by the $E$-Bousfield localization at either the Moore spectrum $\Ss/2$ or the Eilenberg-Maclane spectrum $\HF$. In particular, the 2-completed sphere $(\Ss)^{\smas}_2$ is also an $E_{\infty}$ ring spectrum and admits a good category of ($2$-completed) modules. Denote temporarily its category of modules by $\widehat{\SptC}$. The ring map $\Ss \lto (\Ss)^{\smas}_2$ induces a forgetful functor 
\begin{center}
\begin{tikzpicture}
\matrix (m) [matrix of math nodes, row sep=3em, column sep=4em]
{  \SptC & \widehat{\SptC} \\};
\path[thick, -stealth, font=\small]
(m-1-2.186) edge  (m-1-1);
\end{tikzpicture}
\end{center}
from 2-completed motivic spectra to motivic spectra. As explained in \cite[Section 2.8]{Pelaez}, this forgetful functor creates a symmetric monoidal model structure on $\widehat{\SptC}$.  Moreover, as indicated in the diagram  
\begin{center}
\begin{tikzpicture}
\matrix (m) [matrix of math nodes, row sep=3em, column sep=4em]
{  \SptC & \widehat{\SptC}, \\};
\path[thick, -stealth, font=\tiny]
(m-1-2.186) edge  (m-1-1)
(m-1-1.19) edge node[above] { {\fontsize{5}{5}\selectfont $ - \smas (\Ss)^{\smas}_2 $} }  (m-1-2.168)
(m-1-1.341) edge node[below] { {\fontsize{5}{5}\selectfont $ F\left( (\Ss)^{\smas}_2, - \right) $} }  (m-1-2.203);
\end{tikzpicture}
\end{center}
it is both a left and right Quillen functor via the usual adjunctions. It follows that the forgetful functor preserves all categorical constructions in $\widehat{\SptC}$, i.e., the underlying spectrum of any (co)limit is computed in the underlying category of motivic spectra $\SptC$. We will from now on exclusively work in the 2-completed category without further mention, and drop the completion symbol from the notation. For example, we will denote the category of 2-completed motivic spectra by $\SptC$, the 2-completed motivic sphere spectrum by $\Ss$, the smash product over the 2-completed sphere by $- \smas -$, \ldots etc. With this notation, the motivic Adams spectral sequence produces a non-trivial map $$S^{0,-1} \stackrel{\tau}{\lto} \Ss,$$  
which we can see as being an element in the homotopy groups $\pi_{0,-1}(\Ss)$.


\subsection{The Spectrum $\Ct$ and its Homotopy} \label{sec:Cthtpy}

Recall that we work in a 2-completed setting. Define the 2-cell complex $\Ct$ by the cofiber sequence
\begin{equation} \label{eq:cofseqtau}
S^{0,-1} \stackrel{\tau}{\lto} S^{0,0} \stackrel{i}{\lto} \Ct \stackrel{p}{\lto} S^{1,-1},
\end{equation}
where $i$ denotes the inclusion of its bottom cell and $p$ is the projection on its top cell. Recall from \cite[Section 2.6]{DuggerIsaksenMASS} that the Betti realization functor $\SptC \lto \Spt$ sends the map $\tau$ to the identity $\id$, as shown in the diagram 
\begin{equation*}
\left(S^{0,-1} \stackrel{\tau}{\lto} \Ss \right) \lmapsto \left( S^0 \stackrel{\id}{\lto} S^0 \right).
\end{equation*}
Moreover, it is a left Quillen functor and thus preserves cofiber sequences. This implies that it sends $\Ct$ to a contractible spectrum $\ast \in \Top$ and thus that $\Ct$ is a purely motivic spectrum living in the kernel of Betti realization. Nonetheless, the motivic spectrum $\Ct$ has very tight connections to classical (non-motivic) homotopy theory. Surprisingly, a computation of Hu-Kriz-Ormsby in \cite{HKO11}, allows Isaksen in \cite{StableStems} to express the homotopy groups of this 2-cell complex $\piss(\Ct)$ in terms of the classical Adams-Novikov spectral sequence. Denote by $\Ext^{s,t}_{BP_{\ast}BP}(BP_{\ast},BP_{\ast})$ the $E_2$-page of the classical (2-completed) Adams-Novikov spectral sequence for the topological sphere $S^0$, where as usual $s$ is the Adams filtration and $t$ is the internal degree.

\begin{prop}[{\cite[Proposition 6.2.5]{StableStems}}]
\label{prop:htpycoftau}
The homotopy groups of $\Ct$ are given by
\begin{equation*}
\pi_{s,w}(\Ct) \cong \Ext^{2w-s,2w}_{BP_{\ast}BP}(BP_{\ast},BP_{\ast}) \qquad \qquad \text{for any } s,w \in \Z.
\end{equation*}
\end{prop}

\begin{remark}
Proposition \ref{prop:htpycoftau} is surprising as it is saying that the homotopy groups of a motivic 2-cell complex, which are in principle as complicated to compute as $\piss(\Ss)$, are completely algebraic. More precisely, they are given by the cohomology of the Hopf algebroid $(BP_{\ast}, BP_{\ast}BP)$, which is a very important object in classical chromatic homotopy theory. This bridge allows computations to travel between the classical and the motivic world. See \cite[Chapter 5 and 6]{StableStems} for examples where motivic computations of $\piss(\Ct)$ are used to deduce new information about the classical object $\Ext^{\ast,\ast}_{BP_{\ast}BP}(BP_{\ast},BP_{\ast})$. 
\end{remark}

\begin{remark} \label{rem:CtisEoo}
Since $\Ext^{\ast,\ast}_{BP_{\ast}BP}(BP_{\ast},BP_{\ast})$ admits a natural ring structure, the isomorphism of Proposition \ref{prop:htpycoftau} induces an artificial ring structure on the motivic homotopy groups $\piss(\Ct)$. The starting point of this project was to ask if this induced ring structure of $\piss(\Ct)$ can be realized by a topological ring structure on the spectrum $\Ct$. Even further, the cohomology groups $\Ext^{\ast,\ast}_{BP_{\ast}BP}(BP_{\ast}, BP_{\ast})$ admit higher structure (Massey products, algebraic squaring operations, \ldots) and one can hope that this is the shadow of a highly structured ring multiplication on $\Ct$. We will prove in Section \ref{sec:CtEinfty} that $\Ct$ supports an $E_{\infty}$ ring structure and that the isomorphism 
\begin{equation*}
\piss(\Ct) \cong \Ext^{\ast,\ast}_{BP_{\ast}BP}(BP_{\ast},BP_{\ast})
\end{equation*}
preserves higher products (Toda brackets in homotopy and Massey products in algebra). In other words, the $E_2$-page of the classical Adams-Novikov spectral sequence can be realized with its higher structure as the homotopy of a motivic spectrum.
\end{remark}

The ring structure mentioned in Remark \ref{rem:CtisEoo} will be constructed by obstruction theory. To prepare the computations, we will now deduce some Corollaries about $\piss(\Ct)$ and $\piss( \ECt)$.

\begin{cor}[\cite{GI}] \label{cor:htpycoftau}
The group $\pi_{s,w} ( \Ct )$ is zero when either $w > s$, or $w \leq \frac{1}{2}s$, or $s < 0$, except that $\pi_{0,0} ( \Ct ) \cong \hat{\Z}_2$. This is sketched in Figure \ref{fig:chartCt}.
\begin{proof}
The vanishing regions in $\piss(\Ct)$ come from the vanishing regions of $\Ext^{\ast,\ast}_{BP_{\ast}BP}(BP_{\ast}, BP_{\ast})$ via the isomorphism
\begin{equation*}
\pi_{s,w} (\Ct) \cong \Ext^{2w-s,2w}_{BP_{\ast}BP}(BP_{\ast},BP_{\ast})
\end{equation*}
of Proposition \ref{prop:htpycoftau}. The region $w > s$ corresponds to the vanishing region above the line $t-s =s$ of slope 1 on the $E_2$-page of the Adams-Novikov spectral sequence, the region $w \leq \frac{1}{2}s$ corresponds to the $E_2$-page being 0 in negative Adams filtation $s \leq 0$, and finally $s < 0$ corresponds to $E_2$-page being zero in negative stems $t-s < 0$.
\begin{figure}[ht]
\begin{tikzpicture}[scale=0.2] 
\centering

\draw [ thick, ->, mygray] (0,-4) -- (0,16);
\node [left= 4pt] at (0,16) { $w$ };
\draw [ thick, ->, mygray] (-4,0) -- (25,0);
\node [below=4pt] at (25,0) { $s $ };

\draw [very thick, color3] (0,0) -- (16,16);
\node [right, color3] at (16,15.5) { $w = s$ };
\node [color3] at (20,9) { \textsf{non-vanishing homotopy}  };

\draw [very thick, color3] (0,0) -- (25,7);
\node [right, color3] at (22,5.6) { $w = \frac{1}{2}s$ };

\node at (-2,-2) { \Large \textsf{zero} };
\node at (-2,11) { \Large \textsf{zero} };
\node at (5,11) { \Large \textsf{zero} };
\node at (15,2) { \Large \textsf{zero} };
\node at (15,-2) { \Large \textsf{zero} };
\end{tikzpicture}

\caption{Vanishing regions of the homotopy groups $\pi_{s,w} (\Ct)$.}
\label{fig:chartCt}
\end{figure}
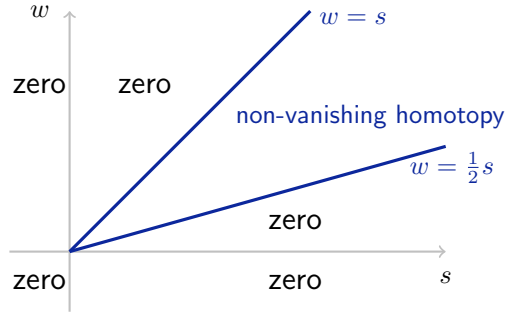
\end{proof}
\end{cor}

\begin{cor} \label{cor:mappingspacecoftau}
The group $\left[ \Sigma^{s,w} \Ct, \Ct \right]$ is zero if either $w > s+2$, or  $w \leq \frac{1}{2}s$, or $s < -1$, except that $\left[ \Ct, \Ct \right] \cong \hat{\Z}_2$ in degree $(0,0)$. This is sketched in Figure \ref{fig:chartmappingCt}.
\begin{proof}
Using the cofiber sequence
\begin{equation*}
S^{s,w} \stackrel{i}{\lto} \Sigma^{s,w} \Ct \stackrel{p}{\lto} S^{s+1,w-1},
\end{equation*}
we get a long exact sequence
\begin{equation*}
\cdots \ltoback \left[ S^{s,w}, \Ct \right] \stackrel{i^{\ast}}{\ltoback} \left[\Sigma^{s,w} \Ct, \Ct \right]  \stackrel{p^{\ast}}{\ltoback} \left[S^{s+1,w-1}, \Ct \right] \ltoback \cdots,
\end{equation*}
after mapping into $\Ct$. The result follows by noticing that the hypothesis of this Corollary force both homotopy groups $\pi_{s,w}\left( \Ct \right)$ and $\pi_{s+1,w-1}\left( \Ct \right)$ to be 0 by the previous Corollary \ref{cor:htpycoftau}.
\begin{figure}[ht]
\begin{tikzpicture}[scale=0.2] 
\centering

\draw [ thick, ->, mygray] (0,-4) -- (0,16);
\node [left= 4pt] at (0,16) { $w$ };
\draw [ thick, ->, mygray] (-6,0) -- (25,0);
\node [below=4pt] at (25,0) { $s $ };

\draw [very thick, color3] (-1,-0.5) -- (-1,1) -- (14,16);
\node [right, color3] at (14,15.5) { $w = s+2$ };
\node [color3] at (16,9) { \textsf{non-vanishing region}  };

\draw [very thick, color3] (-1,-0.5) -- (25,7);
\node [right, color3] at (22,5.6) { $w = \frac{1}{2}s$ };

\node at (-4,-2) { \Large \textsf{zero} };
\node at (-4,11) { \Large \textsf{zero} };
\node at (5,11) { \Large \textsf{zero} };
\node at (15,2) { \Large \textsf{zero} };
\node at (15,-2) { \Large \textsf{zero} };
\end{tikzpicture}

\caption{Vanishing regions of the abelian group $\left[ \Sigma^{s,w} \Ct, \Ct \right]$.}
\label{fig:chartmappingCt}
\end{figure}
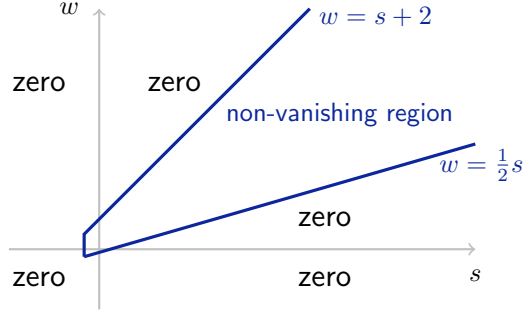

\end{proof}
\end{cor}

\begin{remark}
This result is not sharp and one can slightly improve the non-vanishing region by being careful about choosing which of the 3 conditions of Corollary \ref{cor:htpycoftau} to use. For example, the group $\left[ \Sigma^{-1,0} \Ct, \Ct \right]$ is zero as it sits in a long exact sequence
\begin{equation*}
\cdots \ltoback \pi_{-1,0} \left( \Ct \right) \stackrel{i^{\ast}}{\ltoback} \left[\Sigma^{-1,0} \Ct, \Ct \right]  \stackrel{p^{\ast}}{\ltoback} \pi_{0,-1} \left( \Ct \right) \ltoback \cdots,
\end{equation*}
and both homotopy groups surrounding it are zero. However, none of the 3 conditions of Corollary \ref{cor:mappingspacecoftau} are satisfied for the pair $(s,w) = (-1,0)$ and thus we cannot use it to deduce that $\left[ \Sigma^{-1,0} \Ct, \Ct \right]$ is zero.
\end{remark}

The vanishing of the following groups of homotopy classes of maps will often be used in this document.

\begin{cor} \label{cor:mappingiszero}
The following groups of homotopy classes of maps are zero
\begin{enumerate}
\item $\left[ \Sigma^{0,-1} \Ct, \Ct \right] = 0$,
\item $\left[ \Sigma^{1,0} \Ct, \Ct \right] = 0$,
\item $\left[ \Sigma^{1,-1} \Ct, \Ct \right] = 0$,
\item $\left[ \Sigma^{n,-n} \Ct, \Ct \right] = 0$ for any $n \geq 1$.
\end{enumerate}
\end{cor}

\section{The $E_{\infty}$ Ring Structure on $\Ct$} \label{sec:CtEinfty} 

In this Section we construct the $E_{\infty}$ ring structure on the motivic spectrum $\Ct$. We start by endowing $\Ct$ with a homotopy unital, homotopy associative and homotopy commutative multiplication using elementary techniques with triangulated categories. The $E_{\infty}$ coherences of such a multiplication cannot be constructed by hand via similar techniques and requires some machinery. We will use a version of Robinson's obstruction theory from \cite{RobinsonGamma}, that we adapt to the motivic setting in Section \ref{sec:motoperads}.


\subsection{Motivic $A_{\infty}$ and $E_{\infty}$ Operads and Obstruction Theory} \label{sec:motoperads}

Consider a simplicial symmetric monoidal model category presenting $\SptC$, with smash product $- \smas -$\footnote{For example Jardine's model of motivic symmetric spectra \cite{Jardine}.}, and denote the simplicial mapping space by $\Map(X,Y)$. Given a motivic spectrum $X$, denote its \emph{endomorphism operad} in simplicial sets by $\EEnd(X)$, where $\EEnd(X)_n$ is the simplicial set $\Map(X^{\smas n},X)$. If $F(-,-)$ denotes the internal (motivic) function spectrum, then we recover 
\begin{equation} \label{eq:funcspectrasimp}
\pi_n \left( \EEnd(X)_m \right) \cong \pi_{n,0} \left( F(X^{\smas m},X) \right),
\end{equation}
only exploiting the weight zero homotopy groups of the function spectrum. Fix an $A_{\infty}$ or $E_{\infty}$ operad $\Theta$ in simplicial sets. A \emph{$\Theta$-algebra structure} on a motivic spectrum $X$ is a map of operads
\begin{equation*}
\Theta \lto \EEnd(X).
\end{equation*}
Equivalently, one can see $\Theta$ as an operad in motivic spaces via the constant functor and define a $\Theta$-algebra via the motivic enrichment, which might seem more natural and internal to motivic homotopy theory. Because of this reason, classical (simplicial) operads transported into the motivic world are sometimes called \emph{constant operads}. 

In this paper, we will produce $A_{\infty}$ and $E_{\infty}$ structures by obstruction theory. The obstruction theory for $A_{\infty}$ algebras is well-known, for example \cite[Theorem 3.1]{Angeltveit} (itself inspired by \cite{RobinsonAoo}) exhibits an obstruction class in a certain abelian group. In all our cases, we will show that all the relevant abelian groups for the obstruction theory are zero. The obstruction theory for $E_{\infty}$ algebras is less well-known. We will here briefly recap the work done in \cite{RobinsonGamma} and adapt it to our motivic situation.

We will consider the simplicial $E_{\infty}$ operad $\mathcal{T}$ defined in \cite[Section 5]{RobinsonGamma}. This operad is the product of a combinatorially defined cofibrant simplicial operad with the Barratt-Eccles $E_{\infty}$ (simplicial) operad $E\Sigma_{\bullet}$. It inherits both properties and is thus a cofibrant $E_{\infty}$ operad. The cofibrancy roughly means that the operadic composition maps
\begin{equation} \label{eq:treeopmaps}
\mathcal{T}_n \times \mathcal{T}_m \stackrel{\circ_i}{\ltoo} \mathcal{T}_{m+n-1}
\end{equation}
are injective and that their images intersect in fairly small and regular subcomplexes. We refer to \cite[Section 1.5]{RobinsonWhitehouse} for more details. The injectivity of these maps is a key property that will be used for inductive arguments, since a map out of $\mathcal{T}_{m+n-1}$ is thus already determined on the image of all these composition maps. The bar filtration on the Barrat-Eccles operad induces a filtration on $\mathcal{T}$, where the $n^{\text{th}}$-filtration space of $\mathcal{T}_m$ is denoted by $\mathcal{T}^n_m \subseteq \mathcal{T}_m$. In particular $\mathcal{T}^n_m = \emptyset$ if $n < 0$. Consider now the \emph{diagonal filtration} $\nabla^{\bullet} \mathcal{T}$ which is the sum of the bar filtration from the Barratt-Eccles operad and the filtration by operadic subspaces. More precisely, the $n^{\text{th}}$-graded piece $\nabla^n \mathcal{T} \subset \mathcal{T}$ has $m^{\text{th}}$-space given by $\nabla^n \mathcal{T}_m = \mathcal{T}^{n-m}_m$. If $m >n$, then by definition we have $\nabla^n \mathcal{T}_m = \emptyset$. In particular, observe that $\nabla^n \mathcal{T}$ is not a suboperad as it does not contain $m$-ary operations for $m > n$.

Robinson defines an \emph{$n$-stage for an $E_{\infty}$ structure on $X$} to consist in a map $\nabla^n \mathcal{T} \lto \EEnd(X)$ satisfying some obvious coherences. More precisely, this is the data of $\Sigma_m$-equivariant maps
\begin{equation*}
\mathcal{T}^{n-m}_m \lto \EEnd(X)_m
\end{equation*}
for $0 \leq m \leq n$, which on their restricted domain of definition satisfy the requirements for a morphism of operads. Since the operad $\mathcal{T}$ is non-unital and thus $\mathcal{T}_0 = \mathcal{T}_1 = \emptyset$, we only need to specify these maps for $2 \leq m \leq n$. From the definition of the diagonal filtration one can identify that
\begin{itemize}
\item a 2-stage is the data of a map $\mathcal{T}_2^0 \lto \EEnd(X)_2$, i.e., specifying a map $\mu \colon X \smas X \lto X$,
\item a 3-stage is the data of a 2-stage with the extra structure of an associative and commutative homotopy for the multiplication $\mu$,
\item a 4-stage is the data of a 3-stage with the extra structure of homotopies for the well-known pentagonal and hexagonal axioms \cite{MacLaneAssocComm}, as well as a homotopy saying that the commutativity homotopy itself is homotopy commutative,
\item an $\infty$-stage are the coherences of an $E_{\infty}$ ring structure on $X$ with multiplication $\mu$.
\end{itemize}
An $n$-stage determines an $(n-1)$-stage by restriction, and an $(n-1)$-stage determines an $n$-stage on the boundary $\partial \nabla^n \mathcal{T}$ by injectivity of the composition maps of equation \eqref{eq:treeopmaps}. We refer to \cite[Section 5.2]{RobinsonGamma} for more details. Therefore, given an $(n-1)$-stage, the data of an $n$-stage extending the underlying $(n-1)$-stage consists precisely in the data of extensions
\begin{center}
\begin{tikzpicture}
\matrix (m) [matrix of math nodes, row sep=3em, column sep=4em]
{ \partial \nabla^n \mathcal{T}_m & \nabla^n \mathcal{T}_m \\
                                  & \EEnd(X)_m \\};
\path[thick, -stealth, font=\small]
(m-1-1) edge (m-1-2)
(m-1-1) edge (m-2-2);
\path[thick, dashed, -stealth, font=\small]
(m-1-2) edge (m-2-2);
\end{tikzpicture}
\end{center}
for every $0 \leq m \leq n$. The cofibrancy of the operad $\mathcal{T}$ is used again to show that for any $m$, the map
\begin{equation*}
\partial \nabla^n \mathcal{T}_m \mono \nabla^n \mathcal{T}_m
\end{equation*}
is a principal $\Sigma_m$-equivariant cofibration, whose cofiber is a wedge of spheres $S^{n+2}$ indexed over a set with free $\Sigma_m$-action. 
This allows us to formulate the following result.

\begin{prop} \label{prop:EooObstTheoryMachinery}
Let $X$ be a motivic spectrum with a given $(n-1)$-stage for an $E_{\infty}$ ring structure. 
\begin{enumerate}
\item If the homotopy groups $\pi_{n-3}( \EEnd(X)_m)$ are zero for every $2 \leq m \leq n$, the given $(n-1)$-stage lifts to an $n$-stage.
\item If in addition the homotopy groups $\pi_{n-2}( \EEnd(X)_m)$ are zero for every $2 \leq m \leq n$, the extension is (essentially) unique.
\end{enumerate}
\begin{proof}
The fact that $\partial \nabla^n \mathcal{T}_m \mono \nabla^n \mathcal{T}_m$ is a principal cofibration allows us to rotate it one step to the left, producing the unstable cofiber sequence of simplicial sets
\begin{equation*}
\vee S^{n-3} \lto \partial \nabla^n \mathcal{T}_m \mono \nabla^n \mathcal{T}_m \lto \vee S^{n-2}.
\end{equation*}
An $(n-1)$-stage produces a map $\partial \nabla^n \mathcal{T}_m \lto \EEnd(X)_m$, which extends as in the diagram
\begin{center}
\begin{tikzpicture}
\matrix (m) [matrix of math nodes, row sep=3em, column sep=4em]
{ \vee S^{n-3} & \partial \nabla^n \mathcal{T}_m & \nabla^n \mathcal{T}_m & \vee S^{n-2} \\
                                  & & \EEnd(X)_m & \\};
\path[thick, -stealth, font=\small]
(m-1-1) edge (m-1-2)
(m-1-2) edge (m-1-3)
(m-1-3) edge (m-1-4)
(m-1-2) edge (m-2-3);
\path[thick, dashed, -stealth, font=\small]
(m-1-3) edge (m-2-3);
\end{tikzpicture}
\end{center}
if and only if the relevant composite is zero in the abelian group
\begin{equation*}
\left[ \vee S^{n-3}, \EEnd(X)_m \right] \cong \oplus \pi_{n-3}(\EEnd(X)_m).
\end{equation*}
Moreover, if $\left[S^{n-2}, \EEnd(X)_m \right] = 0$ then the extension is unique up to homotopy.
\end{proof}
\end{prop}

By using equation \eqref{eq:funcspectrasimp} and the fact that a 3-stage is equivalent to a unital, associative and commutative monoid in the homotopy category, we get the following Corollary.

\begin{cor} \label{cor:EooObstTheoryMachinery}
Let $X$ be a motivic spectrum with a map $\mu \colon X \smas X \lto X$ that is homotopy unital, homotopy associative and homotopy commutative.
\begin{enumerate}
\item If the homotopy groups $\pi_{n-3,0}(F(X^{\smas m},X))$ are zero for every $n \geq 4$ and $2 \leq m \leq n$, then $\mu$ can be extended to an $E_{\infty}$ ring structure on $X$.
\item If in addition the homotopy groups $\pi_{n-2,0}( F(X^{\smas m},X))$ are zero for every $n \geq 4$ and $2 \leq m \leq n$, then $\mu$ can be extended to an $E_{\infty}$ ring structure on $X$ in essentially a unique way.
\end{enumerate}
\end{cor}

\begin{remark} \label{rem:RobBetterObst}
These results are extracted from Robinson's work in \cite{RobinsonGamma}, even though they do not explicitly appear in this form in his paper. The reason is because this is not a powerful result when applied to the topological setting for the following reason. Fix a (topological) spectrum $X \in \Spt$. To apply this $E_{\infty}$ obstruction theory to $X$, its endomorphism operad $\EEnd(X)$ has to satisfy the conditions of Proposition \ref{prop:EooObstTheoryMachinery}, which require the homotopy groups $\EEnd(X)_m$ to vanish for all $n \geq 4$ and $2 \leq m \leq n$. In particular, for any fixed $m$ the space $\EEnd(X)_m$ needs to have vanishing homotopy groups in degrees $n \geq m$. The paper \cite{RobinsonGamma} proceeds to study what happens during an extension of an $(n-1)$-stage to an $n$-stage if one allows to perturb underlying stages. This reduces the size of the obstruction groups and gives a constraint between $n$ and $m$, reducing the number of obstruction groups to check. In our motivic setting the obstructions live in the groups $\pi_{n-3,0}( \EEnd(X)_m)$, which are only a small fraction of all homotopy groups $\pisw$. Corollary \ref{cor:EooObstTheoryMachinery} will be sufficient to prove our result.
\end{remark}

\begin{remark}
We should point out that, in analogy with the genuine $G$-equivariant $E_{\infty}$ operads in \cite{BlumbergHill} (called $N_{\infty}$ operads), there ought to be a notion of motivic $A_{\infty}$ and $E_{\infty}$ operads. An algebra over such a motivic operad would have a lot more structure than an algebra over a constant operad, such as transfers upon changing the base scheme. It is possible that such algebras are exactly the objects corresponding to strict commutative ring spectra. However, for the purpose of this paper, constant $A_{\infty}$ and $E_{\infty}$ operads suffice. We will therefore drop the word "constant" and refer to those just as $A_{\infty}$ and $E_{\infty}$ operads.
\end{remark}


\subsection{The Homotopy Ring Structure on $\Ct$} \label{sec:Ctmonoidinhtpycat}

In this Section we construct a ring structure on $\Ct$ up to homotopy. More precisely, we show that $\Ct$ is a unital, associative and commutative monoid in the homotopy category $\HoSptC$. Recall that this is a 3-stage in Robinson's obstruction theory, which can be seen as the initial input to start the obstruction theory. In this Section, we will exclusively work in the stable triangulated category $\HoSptC$, without further mentioning it.

\begin{lemma} \label{lem:multonCt}
There exists a unique left unital multiplication
\begin{equation*}
\Ct \smas \Ct \stackrel{\mu}{\lto} \Ct.
\end{equation*}
\begin{proof}
The equation \eqref{eq:cofseqtau} gives an exact triangle
\begin{equation*}
S^{0,-1} \stackrel{\tau}{\lto} S^{0,0} \stackrel{i}{\lto} \Ct \stackrel{p}{\lto} S^{1,-1},
\end{equation*}
where $i$ denotes the inclusion of the bottom cell and $p$ denotes the projection on the top cell. By smashing it with $ - \smas \Ct$, we get another triangle
\begin{equation*}
S^{0,-1} \smas \Ct \stackrel{\tau}{\lto} S^{0,0} \smas \Ct \stackrel{i_{L}}{\lto} \Ct \smas \Ct \stackrel{p_L}{\lto} S^{1,-1} \smas \Ct,
\end{equation*}
where $i_L$ denotes a left unit and $p_L$ the projection on the top cell of the left factor. Since the abelian group of maps $[\Sigma^{0,-1} \Ct,\Ct] = 0$ by Corollary \ref{cor:mappingiszero}, the map $\tau \in [\Sigma^{0,-1} \Ct,\Ct]$ is zero on $\Ct$. This produces a left unital multiplication $\mu$ on $\Ct$ as shown in the diagram
\begin{equation*}
\begin{tikzpicture}
\matrix (m) [matrix of math nodes, row sep=3em, column sep=4em]
{ S^{0,-1} \wedge \Ct & S^{0,0} \wedge \Ct & \Ct \wedge \Ct & S^{1,-1} \wedge \Ct \\
  & & \Ct. & \\};
\path[thick, -stealth, font=\small]
(m-1-1) edge node[above] {$ \tau $} (m-1-2)
(m-1-2) edge node[above] {$ i_{L} $} (m-1-3)
(m-1-3) edge node[above] {$  p_L $} (m-1-4)
(m-1-2) edge node[auto] {$ \simeq $} (m-2-3);
\path[thick, -stealth, dashed, font=\small]
(m-1-3) edge node[right] {$ \exists \ \mu $} (m-2-3);
\end{tikzpicture}
\end{equation*}
Moreover, since $\left[ \Sigma^{1,-1} \Ct, \Ct \right] = 0$ by Corollary \ref{cor:mappingiszero}, there is no choice for such a map which is unique.
\end{proof}
\end{lemma}

Before studying the properties of this multiplication map $\mu$, we show a fundamental equivalence that will be used throughout the document.

\begin{lemma} \label{lem:equivCtsmasCt}
There is a canonical isomorphism
\begin{equation*}
\Ct \smas \Ct \cong \Ct \vee \SCt.
\end{equation*}
\begin{proof}
Recall that since $[\Sigma^{0,-1} \Ct,\Ct] = 0$, the map $\tau$ is zero on $\Ct$. The exact triangle 
\begin{equation*}
S^{0,-1} \smas \Ct \stackrel{\tau}{\lto} S^{0,0} \smas \Ct \stackrel{i_{L}}{\lto} \Ct \smas \Ct \stackrel{p_L}{\lto} S^{1,-1} \smas \Ct,
\end{equation*}
is thus split, giving both a retraction $\mu$ and a section $s$, as in the diagram
\begin{center}
\begin{tikzpicture}
\matrix (m) [matrix of math nodes, row sep=3em, column sep=4em]
{ S^{0,-1} \wedge \Ct & S^{0,0} \wedge \Ct & \Ct \wedge \Ct & S^{1,-1} \wedge \Ct & \cdots. \\};
\path[thick, -stealth, font=\small]
(m-1-1) edge node[above] {$ \tau = 0 $} (m-1-2)
(m-1-2) edge node[above] {$ i_{L} $} (m-1-3)
(m-1-3) edge node[above] {$  p_L $} (m-1-4)
(m-1-4) edge node[above] {$ \tau = 0 $} (m-1-5);
\path[thick, dotted, -stealth, font=\small]
(m-1-3) edge[bend left] node[below] {$ \exists ! \  \mu $} (m-1-2)
(m-1-4) edge[bend left] node[below] {$ \exists ! \  s $} (m-1-3);
\end{tikzpicture}
\end{center}
As it is the case for $\mu$, the section $s$ is unique since $[\Sigma^{1,-1} \Ct,\Ct] = 0$ by Corollary \ref{cor:mappingiszero}. Moreover, the relation $\mu \comp s \cong 0$ is forced since the composite lives in the zero group $[\Sigma^{1,-1} \Ct,\Ct] = 0$. This gives a canonical identification
\begin{equation*}
\Ct \smas \Ct \cong \Ct \vee \SCt,
\end{equation*}
via the inverse maps
\begin{equation*} 
\Ct \smas \Ct \stackrel{(\mu, p_L)}{\ltooo} \Ct \vee \Sigma^{1,-1} \Ct \qquad \text{ and } \qquad \Ct \vee \Sigma^{1,-1} \Ct \stackrel{i_L + s}{\ltooo} \Ct \smas \Ct.
\end{equation*}
\end{proof}
\end{lemma}

\begin{cor} \label{cor:decCtsmas}
For any $n \geq 2$, there is a canonical isomorphism
\begin{equation*}
\Ct^{\smas n} \cong \bigvee_{i=0}^{n-1} \binom{n-1}{i} \Sigma^{i,-i} \Ct,
\end{equation*}
where we use $\binom{n-1}{i} \Sigma^{i,-i} \Ct$ to indicate a wedge sum of $\binom{n-1}{i}$ terms of the spectrum $\Sigma^{i,-i} \Ct$.
\end{cor}

We will use the identification of Lemma \ref{lem:equivCtsmasCt} to show that $\mu$ endows $\Ct$ with a unital, associative and commutative monoid structure in $\HoSptC$. We first compute the relevant maps on $\Ct \vee \SCt$ after composing with this identification.

\begin{lemma} \label{lemma:idCtsmasCtunitalcommmonid}
After the canonical identification $\Ct \smas \Ct \cong \Ct \vee \SCt$ of Lemma \ref{lem:equivCtsmasCt}
\begin{enumerate}
\item the multiplication map $\Ct \smas \Ct \stackrel{\mu}{\lto} \Ct$ is given by the matrix $$\Ct \vee \SCt \stackrel{\left[ \begin{smallmatrix} \id&0 \end{smallmatrix} \right]}{\ltoo} \Ct,$$ i.e., by the canonical projection onto the first factor,

\item the factor swap map $\Ct \smas \Ct \stackrel{\chi}{\lto} \Ct \smas \Ct$ is given by the matrix $$\Ct \vee \Sigma^{1,-1} \Ct \stackrel{\left[ \begin{smallmatrix} \id&0\\ i \comp p&-\id \end{smallmatrix} \right]}{\ltooo} \Ct \vee \Sigma^{1,-1} \Ct.$$
\end{enumerate}
\begin{proof}
\mbox{}
\begin{enumerate}
\item The composite
\begin{equation*}
\Ct \vee \Sigma^{1,-1} \Ct \stackrel{i_L + s}{\ltooo} \Ct \smas \Ct \stackrel{\mu}{\lto} \Ct
\end{equation*}
restricts to the identity on $\Ct$ since $\mu$ is a retraction of $i_L$, and to zero on $\SCt$ since $s \comp \mu = 0$ by Lemma \ref{lem:equivCtsmasCt}.

\item We claim that the following diagram 
\begin{center}
\begin{tikzpicture} [ampersand replacement=\&]
\matrix (m) [matrix of math nodes, row sep=3em, column sep=4em]
{  \CtCt \& \CtCt \\
   \Ct \vee \Sigma^{1,-1} \Ct \& \Ct \vee \Sigma^{1,-1} \Ct \\};
\path[thick, -stealth]
(m-1-1) edge node[above] {$ \chi $} (m-1-2)
(m-2-1) edge node[left] {$ i_L + s $} (m-1-1)
(m-1-2) edge node[right] {$ ( \mu, p_L) $} (m-2-2);
\path[thick, -stealth, dashed]
(m-2-1) edge node[above] {$ \left[ \begin{smallmatrix} \id&0\\ i \comp p&-\id \end{smallmatrix}  \right] $} (m-2-2);
\end{tikzpicture}
\end{center}
commutes. First observe that the top right entry is forced to be zero since $\left[ \Sigma^{1,-1} \Ct, \Ct \right] = 0$ by Corollary \ref{cor:mappingiszero}. The bottom left entry can be computed explicitly by a simple diagram chase. It is
\begin{equation*}
\Ss \smas \Ct \stackrel{i \smas \id}{\ltooo} \Ct \smas \Ct \stackrel{\chi}{\ltooo} \Ct \smas \Ct \stackrel{p \smas \id}{\ltooo} S^{1,-1} \smas \Ct,
\end{equation*}
which is homotopic to the composite
\begin{equation*}
\Ss \smas \Ct \stackrel{\chi}{\lto} \Ct \smas \Ss \stackrel{\id \smas i}{\ltooo} \Ct \smas \Ct \stackrel{p \smas \id}{\ltooo} S^{1,-1} \smas \Ct.
\end{equation*}
By commuting $\id \smas i$ and $p \smas \id$ and using the canonical equivalences $\Ss \smas \Ct = \Ct = \Ct \smas \Ss$ we can rewrite it as
\begin{equation*}
\Ct \stackrel{p}{\lto} S^{1,-1} \stackrel{i}{\lto} \Sigma^{1,-1} \Ct.
\end{equation*}
For the diagonal entries, recall that $\left[ \Ct, \Ct \right] \cong \hat{\Z}_2$ and that the matrix has to be an involution since $\chi$ is. This forces the diagonal entries to be $+\id$ and $-\id$. One could conclude by arguing that the top left entry arises by commuting $\Ct$ with $\Ss$, and thus should be $+\id$, while the bottom right entry arises by commuting $\Ct$ with $S^{1,-1}$, and thus should be $-\id$. More precisely, consider the diagram
\begin{center}
\begin{tikzpicture}
\matrix (m) [matrix of math nodes, row sep=3em, column sep=3em] 
{ \Ss \smas \Ss & \Ct \smas \Ct \\
   \Ss & \Ct. \\};
\path[thick, -stealth]
(m-1-1) edge node[above] {$ i \smas i $} (m-1-2)
(m-1-1) edge node[left] {$ \cong $} (m-2-1)
(m-1-2) edge node[right] {$ \mu $} (m-2-2)
(m-2-1) edge node[above] {$ i $} (m-2-2);
\end{tikzpicture}
\end{center}
By factoring the map $i \smas i$ as $\id \smas i$ followed by $i_L = i \smas \id$, and using that $\mu \comp i_L = \id$, one sees that the diagram commutes up to the usual canonical equivalences of smashing with $\Ss$. By factoring it the other way now, as $i \smas \id$ followed by $\id \smas i$, we get that $\mu \comp (\id \smas i) = \id$. This shows that the top left entry of the matrix is $\id$. The bottom right entry is thus forced to be $-\id$ since the matrix is an involution.
\end{enumerate}
\end{proof}
\end{lemma}

\begin{prop} \label{prop:Ctismonoidinhtpy}
The unique left unital multiplication map $\Ct \smas \Ct \stackrel{\mu}{\lto} \Ct$ turns $\Ct$ into a unital, associative and commutative monoid in $\HoSptC$.
\begin{proof}
Consider the diagram 
\begin{equation} \label{diag:Ctisunitalcomm}
\begin{tikzpicture}
\matrix (m) [matrix of math nodes, row sep=2em, column sep=3em] 
{ & \Ct \smas \Ct & \Ct \smas \Ct &  \\
  \Ct &  &  & \Ct, \\
  & \Ct \vee \SCt & \Ct \vee \SCt & \\};
\path[thick, -stealth]
(m-2-1) edge node[auto] {$ i_L $} (m-1-2)
(m-1-2) edge node[above] {$ \chi $} (m-1-3)
(m-1-3) edge node[auto] {$ \mu $} (m-2-4)
(m-1-2) edge node[right = 25pt, above] {$ (\mu,p\smas \id) $} (m-3-2)
(m-3-3) edge node[left = 18pt, above] {$ i_L + s $} (m-1-3)
(m-3-2) edge node[above = 3pt] {$ \left[ \begin{smallmatrix} \id&0\\ i \comp p&-\id \end{smallmatrix}  \right] $} (m-3-3)
(m-3-3) edge node[right = 10pt, below] {$ \left[ \begin{smallmatrix} \id&0 \end{smallmatrix}  \right] $} (m-2-4);
\path[thick,dashed, -stealth]
(m-2-1) edge (m-3-2);
\end{tikzpicture}
\end{equation}
which is commutative by Lemma \ref{lemma:idCtsmasCtunitalcommmonid}. Since $\mu$ is left unital and since $p \comp i = 0$, the dashed arrow is given by the canonical inclusion. It follows that the composite $\mu \comp \chi \comp i_L$ is simply given by the matrix multiplication 
\begin{equation*}
\left[ \begin{smallmatrix} \id &0 \end{smallmatrix} \right] \cdot \left[ \begin{smallmatrix} \id&0\\ i \comp p &-\id \end{smallmatrix} \right] \cdot \left[ \begin{smallmatrix} \id\\ 0 \end{smallmatrix} \right] = \id.
\end{equation*}
Since the right unit is given by $\chi \comp i_L$, this shows that $\mu$ is right unital. To show that $\mu$ is commutative, we have to compute the composite 
\begin{equation*}
\Ct \smas \Ct \stackrel{\chi}{\lto} \Ct \smas \Ct \stackrel{\mu}{\lto} \Ct.
\end{equation*}
We can again read it from diagram \eqref{diag:Ctisunitalcomm}, where it is given by the matrix multiplication
\begin{equation*}
\left[ \begin{smallmatrix} \id &0 \end{smallmatrix} \right] \cdot \left[ \begin{smallmatrix} \id&0\\ i \comp p &-\id \end{smallmatrix} \right] \cdot \left[ \begin{smallmatrix} \mu \\ p \smas \id \end{smallmatrix} \right] = \mu,
\end{equation*}
showing that $\mu$ is commutative. To see that $\mu$ is associative, we will show that the map
\begin{equation*}
\Ct \smas \Ct \smas \Ct \stackrel{\mu \comp (1 \smas \mu - \mu \smas 1)}{\ltoooo} \Ct
\end{equation*}
is zero. By left and right unitatlity it restricts to zero on the subspectrum
\begin{equation} \label{eq:ressubspectrum}
\left( \Ss \smas \Ct \smas \Ct \right) \vee \left( \Ct \smas \Ss \smas \Ct \right) \vee \left( \Ct \smas \Ct \smas \Ss \right) \inj \Ct \smas \Ct \smas \Ct.
\end{equation}
By \cite[Lemma 3.6]{Strickland}, there is a bijection between maps $\Ct \smas \Ct \smas \Ct \lto \Ct$ that restrict to zero on the subspectrum of equation \eqref{eq:ressubspectrum}, and maps
\begin{equation*}
S^{3,-3} = S^{1,-1} \smas S^{1,-1} \smas S^{1,-1} \lto \Ct.
\end{equation*}
Here $S^{1,-1}$ appears because it is the cofiber of the unit map $\Ss \lto \Ct$. By Corollary \ref{cor:htpycoftau}, we have that $\pi_{3,-3}(\Ct) = 0$, which shows that there is a unique such map. Since the zero map $\Ct \smas \Ct \smas \Ct \lto \Ct$ restricts to zero on the subspectrum of equation \eqref{eq:ressubspectrum}, it is the unique such map. This shows that $\mu \comp (1 \smas \mu - \mu \smas 1)$ is zero, i.e., that $\mu$ is associative.
\end{proof}
\end{prop}


\subsection{The $E_{\infty}$ Ring Structure on $\Ct$}

In this Section, we will use Robinson's obstruction theory from Section \ref{sec:motoperads} to construct the $E_{\infty}$ ring structure on $\Ct$. In the previous Section \ref{sec:Ctmonoidinhtpycat} we endowed $\Ct$ with a unital, associative and commutative monoid structure in the the homotopy category $\HoSptC$. Recall that this to a 3-stage in Robinson's obstruction theory. We will now use Corollary \ref{cor:EooObstTheoryMachinery} to rigidify this multiplication to an $E_{\infty}$ ring structure in $\SptC$. Although not needed for the $E_{\infty}$ ring structure, as a warm-up, we first show in Proposition \ref{prop:multhtpyassoc} that $\Ct$ admits a unique $A_{\infty}$ ring structure.

\begin{prop} \label{prop:multhtpyassoc}
The multiplication $\mu$ on $\Ct$ can be uniquely extended to an $A_{\infty}$ multiplication.
\begin{proof}
An $A_2$ structure corresponds to unital homotopies (left and right), and an $A_3$ structure adds an associative homotopy. We constructed both structures in Proposition \ref{prop:Ctismonoidinhtpy}. The $A_{\infty}$ obstruction theory originated in \cite{RobinsonAoo} exhibits obstruction classes to extend an $A_{n-1}$ structure to an $A_n$ structure. In more modern language, \cite[Theorem 3.1]{Angeltveit} exhibits the obstruction to go from $A_{n-1}$ structure to an $A_n$ structure as an element in the abelian group
\begin{equation} \label{eq:Angeltveitsgroup}
\left[ \Sigma^{n-3,0} {S^{n,-n}}, \Ct \right] \cong \left[ S^{2n-3, -n}, \Ct \right] = \pi_{2n-3,-n} (\Ct).
\end{equation}
Corollary \ref{cor:htpycoftau} shows that these groups are zero for any $n$ (we really just need $n \geq 4$), which shows that $\mu$ can be extended to an $A_{\infty}$ structure. Furthermore, given that an $A_{n-1}$ structure extends to $A_n$ structure, the possible extensions are in bijection with the abelian group
\begin{equation*}
\left[ \Sigma^{n-2,0} {S^{n,-n}}, \Ct \right] \cong \pi_{2n-2,-n} (\Ct).
\end{equation*}
This group is also zero for any $n$, showing that $\mu$ can be uniquely extended to an $A_{\infty}$ structure.
\end{proof}
\end{prop}

\begin{remark}
For the case $n =3$, i.e., to endow $\Ct$ with an $A_3$ structure, the obstruction group from equation \eqref{eq:Angeltveitsgroup} is $\pi_{3,-3}(\Ct)$. Observe that this is the exact same group that appears in Proposition \ref{prop:Ctismonoidinhtpy}, where we show with elementary techniques that $\Ct$ admits an $A_3$ structure.
\end{remark}

\begin{remark} \label{rem:Mahowaldconj}
Mahowald conjectured that no non-trivial topological 2-cell complex posses an $A_{\infty}$ structure. There are 2 trivial cases to exclude which are the cofiber of the zero map and the cofiber of the identity map, as shown in the cofiber sequences
\begin{equation*}
S^0 \stackrel{0}{\lto} S^0 \lto S^1 \vee S^0 \qquad \text{ and } \qquad S^0 \stackrel{\id}{\lto} S^0 \lto \ast.
\end{equation*}
Since motivic spheres Betti realize to topological spheres, motivic 2-cell complexes Betti realize to topological 2-cell complexes. Moreover, since we are using simplicial (constant) operads, motivic algebras over $A_n$ or $E_n$ operads realize to classical algebras over the same $A_n$ or $E_n$ operads. However, the fact that $\Ct$ admits an $A_{\infty}$ ring structure does not contradict Mahowald's conjecture, as the map $S^{0,-1} \stackrel{\tau}{\lto} \Ss$ realizes to the identity map $S^0 \stackrel{\id}{\lto} S^0$.
\end{remark}

\begin{thm} \label{thm:multCthtpycom}
The multiplication $\mu$ on $\Ct$ can be uniquely extended to an $E_{\infty}$ multiplication.
\begin{proof}
We showed in Proposition \ref{prop:Ctismonoidinhtpy} that $\Ct$ is a unital, associative and commutative monoid in the homotopy category $\HoSptC$. This corresponds to a 3-stage in Robinson's obstruction theory. By Corollary \ref{cor:EooObstTheoryMachinery}, the obstructions of extending this 3-stage to an $E_{\infty}$ ring structure live in
\begin{equation*}
\pi_{n-3,0}( F(\Ct^{\smas m}, \Ct) ) \cong \left[ \Sigma^{n-3,0} \Ct^{\smas m}, \Ct \right]
\end{equation*}
for $n \geq 4$ and $2 \leq m \leq n$. Recall from Corollary \ref{cor:mappingspacecoftau} that $\left[ \Sigma^{s,w} \Ct, \Ct \right]$ has in particular a vanishing region for $s \geq 0$ and $2w \leq s$. We now show that all obstruction groups live in this vanishing area. By the equivalence 
\begin{equation*} 
\left[ \Ct^{\smas m}, \Ct \right] \cong \bigoplus_{i=0}^{m-1} \binom{m-1}{i} \left[ \Sigma^{i,-i} \Ct, \Ct \right]
\end{equation*}
of Corollary \ref{cor:decCtsmas}, we have
\begin{equation*}
\pi_{n-3,0}( F(\Ct^{\smas m}, \Ct) ) \cong \left[ \Sigma^{n-3,0} \Ct^{\smas m}, \Ct \right] \cong \bigoplus_{i=0}^{m-1} \binom{m-1}{i} \left[ \Sigma^{n-3 +i,-i} \Ct, \Ct \right].
\end{equation*}
In particular, all the obstructions live in groups of the form $\left[ \Sigma^{s,w} \Ct, \Ct \right]$ where the $s$-coordinate satisfies
\begin{equation*}
s = n - 3 + i \geq  4 - 3 + i \geq 1 
\end{equation*}
while the $w$-coordinate satisfies both
\begin{equation*}
w = -i \leq 0 \qquad \text{ and } \qquad w = -i = n- s -3 \geq 1 -s.
\end{equation*}
This corresponds to the region bounded by $s \geq 1$ and  $1 -s \leq w \leq s$, which lies entirely in the vanishing area described above. The situation is summarized in Figure \ref{fig:mappingspaceCtobstructionsforEoo}. Similarly, recall from Corollary \ref{cor:EooObstTheoryMachinery} that the obstructions for uniqueness of such an $E_{\infty}$ ring structure live in groups of the form
\begin{equation*}
\pi_{n-2,0}( F(\Ct^{\smas m}, \Ct) ) \cong \left[ \Sigma^{n-2,0} \Ct^{\smas m}, \Ct \right].
\end{equation*}
A similar analysis shows that all obstruction groups again live in the vanishing region, as described in Figure \ref{fig:mappingspaceCtobstructionsforEoo}. This shows that $\Ct$ admits a unique $E_{\infty}$ ring structure.
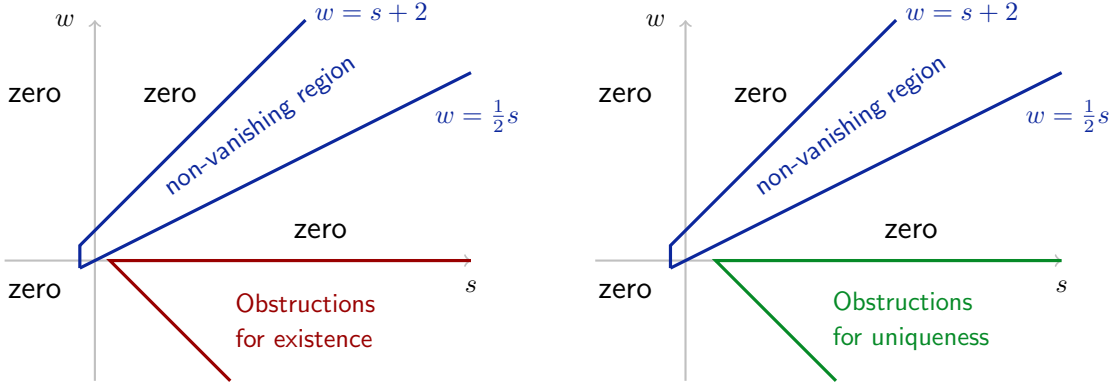
\begin{figure}[ht]
\begin{tikzpicture}[scale=0.2] 
\centering

\draw [ thick, ->, mygray] (0,-8) -- (0,16);
\node [left= 4pt] at (0,16) { $w$ };
\draw [ thick, ->, mygray] (-6,0) -- (25,0);
\node [below=4pt] at (25,0) { $s $ };

\draw [very thick, color3] (-1,-0.5) -- (-1,1) -- (14,16);
\node [right, color3] at (14,16.5) { $w = s+2$ };
\node [color3, rotate=35] at (11,9) { \textsf{non-vanishing region}  };

\draw [very thick, color3] (-1,-0.5) -- (25,12.5);
\node [right, color3]  at (22,9.6) { $w = \frac{1}{2}s$ };

\node at (-4,-2) { \Large \textsf{zero} };
\node at (-4,11) { \Large \textsf{zero} };
\node at (5,11) { \Large \textsf{zero} };
\node at (15,2) { \Large \textsf{zero} };

\draw [very thick, color4] (25,0) -- (1,0) -- (9,-8);
\node [color4, align=left] at (14,-4) { \textsf{Obstructions} \\ \textsf{for existence}  };
\end{tikzpicture}
\qquad
\begin{tikzpicture}[scale=0.2] 
\centering

\draw [ thick, ->, mygray] (0,-8) -- (0,16);
\node [left= 4pt] at (0,16) { $w$ };
\draw [ thick, ->, mygray] (-6,0) -- (25,0);
\node [below=4pt] at (25,0) { $s $ };

\draw [very thick, color3] (-1,-0.5) -- (-1,1) -- (14,16);
\node [right, color3] at (14,16.5) { $w = s+2$ };
\node [color3, rotate=35] at (11,9) { \textsf{non-vanishing region}  };

\draw [very thick, color3] (-1,-0.5) -- (25,12.5);
\node [right, color3]  at (22,9.6) { $w = \frac{1}{2}s$ };

\node at (-4,-2) { \Large \textsf{zero} };
\node at (-4,11) { \Large \textsf{zero} };
\node at (5,11) { \Large \textsf{zero} };
\node at (15,2) { \Large \textsf{zero} };

\draw [very thick, color2] (25,0) -- (2,0) -- (10,-8);
\node [color2, align=left] at (15,-4) { \textsf{Obstructions} \\ \textsf{for uniqueness}  };
\end{tikzpicture}
\caption{Chart of $\left[ \Sigma^{s,w} \Ct, \Ct \right]$ where all obstruction groups live in the vanishing region.}  \label{fig:mappingspaceCtobstructionsforEoo}
\end{figure}
\end{proof}
\end{thm}

\begin{cor} \label{cor:highlystructisoCtAN}
There is an isomorphism of rings
\begin{equation*}
\piss(\Ct) \cong \Ext^{\ast,\ast}_{BP_{\ast}BP}(BP_{\ast},BP_{\ast}),
\end{equation*}
which sends Massey products in $\Ext$ to Toda brackets in $\piss$, and vice-versa.
\begin{proof}
Since $\Ct$ is an $E_{\infty}$ ring spectrum, its motivic Adams-Novikov spectral sequence is multiplicative and converges to an associated graded of the ring $\pi_{\ast,\ast} (\Ct)$. Recall from Proposition \ref{prop:htpycoftau} that the spectral sequence collapses at $E_2$ with no possible hidden extensions as a module over the spectral sequence for $\Ss$. For the exact same reason, there are no possible hidden extensions as a multiplicative spectral sequence. By the Moss convergence theorem \cite{Moss}, we get a highly structured bigraded isomorphism 
\begin{equation} \label{eq:CtBPeq1}
\Ext_{BPGL_{\ast,\ast}BPGL}(BPGL_{\ast,\ast},BPGL_{\ast,\ast}/\tau) \cong \piss(\Ct),
\end{equation}
between the $E_2$-page and the output of the spectral sequence. More precisely, Massey products computed in $\Ext$ converge to Toda brackets computed in $\piss(\Ct)$.

Until the end of the proof, denote the motivic Brown-Peterson spectrum $BPGL$ by $B$. To finish the proof, we have to show that there is a highly structured ring isomorphism
\begin{equation*}
\Ext_{B_{\ast,\ast}B/\tau}(B_{\ast,\ast}/\tau,B_{\ast,\ast}/\tau) \cong \Ext_{B_{\ast,\ast}B}(B_{\ast,\ast},B_{\ast,\ast}/\tau).
\end{equation*}
These are $\Ext$-groups computed in comodules and since the first variable is projective (even free) over the base ring, both of those $\Ext$ terms can be computed from their cobar complex \cite[Corollary A1.2.12]{Ravenel}. Moreover, since the cobar complex also controls the Massey products in the $\Ext$-ring, this will give an isomorphism preserving this structure. The cobar complex of the left $\Ext$-group is given by
\begin{equation*}
{B_{\ast,\ast}/\tau} \otimes_{B_{\ast,\ast}/\tau} {B_{\ast,\ast}B/\tau} \otimes_{B_{\ast,\ast}/\tau} {B_{\ast,\ast}/\tau} \lto {B_{\ast,\ast}/\tau} \otimes_{B_{\ast,\ast}/\tau} {B_{\ast,\ast}B/\tau} \otimes_{B_{\ast,\ast}/\tau} {B_{\ast,\ast}B/\tau} \otimes_{B_{\ast,\ast}/\tau} {B_{\ast,\ast}/\tau} \lto \cdots,
\end{equation*}
while the cobar complex of the right term is given by
\begin{equation*}
{B_{\ast,\ast}} \otimes_{B_{\ast,\ast}} {B_{\ast,\ast}B} \otimes_{B_{\ast,\ast}} {B_{\ast,\ast}/\tau} \lto {B_{\ast,\ast}} \otimes_{B_{\ast,\ast}} {B_{\ast,\ast}B} \otimes_{B_{\ast,\ast}} {B_{\ast,\ast}B} \otimes_{B_{\ast,\ast}} {B_{\ast,\ast}/\tau} \lto \cdots.
\end{equation*}
By iterating the ring isomorphism
\begin{equation*}
{B_{\ast,\ast}/\tau} \otimes_{B_{\ast,\ast}/\tau} {B_{\ast,\ast}B/\tau} \cong B_{\ast,\ast} \otimes_{B_{\ast,\ast}} {B_{\ast,\ast}/\tau},
\end{equation*}
these cobar complexes are isomorphic as dga's. By taking cohomology, we get an isomorphism
\begin{equation} \label{eq:CtBPeq2}
\Ext_{B_{\ast,\ast}B}(B_{\ast,\ast},B_{\ast,\ast}/\tau) \cong \Ext_{B_{\ast,\ast}B/\tau}(B_{\ast,\ast}/\tau,B_{\ast,\ast}/\tau)
\end{equation}
that preserves Massey products. The trigraded $\Ext$-term $\Ext_{B_{\ast,\ast}B/\tau}(B_{\ast,\ast}/\tau,B_{\ast,\ast}/\tau)$ is really bigraded because of the relation $t = 2w$ between the internal degree $t$ and the weight $w$. Therefore, when working mod $\tau$, we can regrade everything in sight by keeping the internal degree and forgetting the weight. With this convention, the degree of $v_n \in B_{\ast}/\tau$ is the single number $2^{n+1} -2$ and thus there is an isomorphism of Hopf algebroids $B_{\ast}B/\tau \cong BP_{\ast}BP$. This provides the (higher) ring isomorphism
\begin{equation} \label{eq:CtBPeq3}
\Ext_{B_{\ast}B/\tau}(B_{\ast}/\tau,B_{\ast}/\tau) \cong \Ext_{BP_{\ast}BP}(BP_{\ast},BP_{\ast}).
\end{equation}
By combining the isomorphisms of equation \eqref{eq:CtBPeq1}, \eqref{eq:CtBPeq2} and \eqref{eq:CtBPeq3}, we get an isomorphism
\begin{equation*}
\piss(\Ct) \cong \Ext_{BP_{\ast}BP}(BP_{\ast},BP_{\ast})
\end{equation*}
of higher rings, that sends Toda brackets to Massey products and vice-versa.
\end{proof}
\end{cor}

\section{(Co-)operations on $\Ct$ } \label{sec:computationsCt}

In this Section we describe the homotopy types of $\Ct \smas \Ct$ and $\ECt$ as ring spectra. Understanding their homotopy types is crucial for the computation of the Steenrod algebra of the spectrum $\HF \smas \Ct$ in Section \ref{sec:HFtwomodtau}. Most proofs are done by diagram chasing and identifying composites of maps.

\subsection{The Spectrum $\Ct \smas \Ct$}

The $E_{\infty}$ ring structure on $\Ct$ induces an $E_{\infty}$ ring structure on the smash product $\Ct \smas \Ct$ via the multiplication
\begin{equation*}
\mu_{\Ct \smas \Ct} \colon (\Ct \smas \Ct) \smas (\Ct \smas \Ct) \stackrel{1 \smas \chi \smas 1}{\ltooo} \Ct \smas \Ct \smas \Ct \smas \Ct \stackrel{\mu \smas \mu}{\ltoo} \Ct \smas \Ct.
\end{equation*}
Here $\mu$ denotes the multiplication map on $\Ct$ and $\chi$ denotes the factor swap map. Recall from Lemma \ref{lem:equivCtsmasCt} that there is a canonical equivalence
\begin{equation*}
\Ct \smas \Ct \simeq \Ct \vee \Sigma^{1,-1} \Ct,
\end{equation*}
describing the additive homotopy type of $\Ct \smas \Ct$. The next lemma describes its ring structure.

\begin{lemma} \label{lem:multCtsmasCt}
Under the canonical vertical identifications given by
\begin{center}
\begin{tikzpicture}
\matrix (m) [matrix of math nodes, row sep=2em, column sep=4em]
{  (\Ct \smas \Ct)  \smas (\Ct \smas \Ct) & \Ct \smas \Ct  \\
   (\Ct \vee \SCt) \smas (\Ct \vee \SCt) & \Ct \vee \SCt \\
    (\Ct \smas \Ct)  \vee (\SCt \smas \Ct) \vee  (\Ct \smas \SCt) \vee (\SCt \smas \SCt) &  \Ct \vee \SCt,  \\};
\path[thick, -stealth]
(m-1-1) edge node[above] {$ \mu_{\Ct \smas \Ct} $} (m-1-2)
(m-2-1) edge (m-2-2);
\path[thick, -stealth, dashed]
(m-3-1) edge (m-3-2);
\path[thick]
(m-1-1) edge node[left] {$ \simeq $} (m-2-1)
(m-2-2) edge node[right] {$ \simeq $} (m-1-2)
(m-2-1) edge node[left] {$ = $} (m-3-1)
(m-2-2) edge node[right] {$ = $} (m-3-2);
\end{tikzpicture}
\end{center}
the multiplication on $\Ct \smas \Ct$ is given by the maps
\begin{align*}
\Ct \smas \Ct &\stackrel{(\mu,0)}{\ltooo} \Ct \vee \SCt \\
\SCt \smas \Ct &\stackrel{(0,\mu)}{\ltooo} \Ct \vee \SCt \\
\Ct \smas \SCt &\stackrel{(0, \mu)}{\ltooo} \Ct \vee \SCt \\
\SCt \smas \SCt &\stackrel{(0,0)}{\ltooo} \Ct \vee \SCt.
\end{align*}
\begin{proof}
These four maps are given by a simple diagram chase, where we only have to be careful with the identifications. For simplicity, let's denote the sphere spectrum $\Ss$ by $S$, and ignore or denote by $1$ some identity maps $\id$ in the following diagrams. Recall the cofiber sequence
\begin{equation*}
S^{0,-1} \stackrel{\tau}{\lto} \Ss \stackrel{i}{\lto} \Ct \stackrel{p}{\lto} S^{1,-1}
\end{equation*}
from equation \eqref{eq:cofseqtau}. The first map $\Ct \smas \Ct \lto \Ct \vee \SCt$ corresponds to the composite 
\begin{equation*}
(\mu, p\smas 1) \comp (\mu \smas \mu) \comp (1 \smas \chi \smas 1) \comp (i \smas i),
\end{equation*}
which is embedded in the commutative diagram
\begin{center}
\begin{tikzpicture}
\matrix (m) [matrix of math nodes, row sep=2em, column sep=4em]
{  (\Ct \smas \Ct)  \smas (\Ct \smas \Ct) & \Ct \smas \Ct \smas \Ct \smas \Ct & \Ct  \smas \Ct & \Ct \vee \SCt  \\
   (S \smas \Ct) \smas (S \smas \Ct) &  S \smas S \smas \Ct \smas \Ct & S \smas \Ct \smas \Ct. & \\ };
\path[thick, -stealth]
(m-1-1) edge node[above] {$ 1 \smas \chi \smas 1 $} (m-1-2)
(m-1-2) edge node[above] {$ \mu \smas \mu $} (m-1-3)
(m-1-3) edge node[above] {$ (\mu, p\smas 1) $} (m-1-4)
(m-2-1) edge node[auto] {$ i \smas i $} (m-1-1)
(m-2-2) edge node[auto] {$ i \smas i $} (m-1-2)
(m-2-3) edge node[auto] {$ i  \smas \mu $} (m-1-3);
\path[thick]
(m-2-1) edge node[above] {$ \simeq $} (m-2-2)
(m-2-2) edge node[above] {$ \simeq $} (m-2-3);
\end{tikzpicture}
\end{center}
We can compute by the other path, where we use that the map
\begin{equation*}
S \smas \Ct \smas \Ct \stackrel{i \smas \mu}{\ltooo} \Ct \smas \Ct
\end{equation*}
decomposes as
\begin{equation*}
S \smas \Ct \smas \Ct \stackrel{1 \smas \mu}{\ltooo} S \smas \Ct \stackrel{i \smas 1}{\ltoo} \Ct \smas \Ct,
\end{equation*}
and by using that $p \comp  i = 0$ and $\mu \comp (i \smas 1) = \id$. For the second map, the canonical splitting of Lemma \ref{lem:equivCtsmasCt} induces a splitting
\begin{equation*}
\SCt \smas \Ct \simeq \SCt \vee \SSCt.
\end{equation*}
By Corollary \ref{cor:mappingiszero} we have $\left[ \SCt, \Ct \right] = \left[ \SSCt, \Ct \right] = 0$, and thus the second map
\begin{equation*}
\SCt \smas \Ct \lto \Ct \vee \SCt
\end{equation*}
corestricts to zero on $\Ct$. To compute the other part, recall first from Lemma \ref{lem:equivCtsmasCt} that the map $p \smas 1$ admits a canonical section $s$, as shown in the cofiber sequence
\begin{center}
\begin{tikzpicture}
\matrix (m) [matrix of math nodes, row sep=3em, column sep=4em]
{ S^{0,-1} \wedge \Ct & S \wedge \Ct & \Ct \wedge \Ct & S^{1,-1} \wedge \Ct & \cdots. \\};
\path[thick, -stealth, font=\small]
(m-1-1) edge node[above] {$ \tau = 0 $} (m-1-2)
(m-1-2) edge node[above] {$ i_{L} $} (m-1-3)
(m-1-3) edge node[above] {$  p_L $} (m-1-4)
(m-1-4) edge node[above] {$ \tau = 0 $} (m-1-5);
\path[thick, dotted, -stealth, font=\small]
(m-1-3) edge[bend left] node[below] {$ \exists ! \  \mu $} (m-1-2)
(m-1-4) edge[bend left] node[below] {$ \exists ! \  s $} (m-1-3);
\end{tikzpicture}
\end{center}
The second map is the composite in the commutative diagram
\begin{center}
\begin{tikzpicture}
\matrix (m) [matrix of math nodes, row sep=2em, column sep=4em]
{  (\Ct \smas \Ct)  \smas (\Ct \smas \Ct) & \Ct \smas \Ct \smas \Ct \smas \Ct & \Ct  \smas \Ct & \SCt  \\
   (S^{1,-1} \smas \Ct) \smas (S \smas \Ct) &  \Ct \smas \Ct \smas S \smas \Ct & \Ct \smas S \smas \Ct \smas \Ct. &  \\ };
\path[thick, -stealth]
(m-1-1) edge node[above] {$ 1 \smas \chi \smas 1 $} (m-1-2)
(m-1-2) edge node[above] {$ \mu \smas \mu $} (m-1-3)
(m-1-3) edge node[above] {$ p \smas 1 $} (m-1-4)
(m-2-1) edge node[auto] {$ s \smas (i \smas 1) $} (m-1-1)
(m-2-3) edge node[right] {$ 1 \smas \mu $} (m-1-3)
(m-2-1) edge node[below] {$ s \smas (1 \smas 1) $} (m-2-2);
\path[thick]
(m-2-2) edge node[above] {$ \simeq $} (m-2-3);
\path[thick, -stealth, dashed]
(m-2-2) edge node[right = 3 pt, above = 0 pt] {$  i  $} (m-1-1)
(m-2-3) edge node[right = 10 pt, above = 0pt] {$  i  $} (m-1-2);
\end{tikzpicture}
\end{center}
We again compute it by following the other path
\begin{equation*}
(p \smas 1) \comp (1 \smas \mu) \comp (s \smas (1 \smas 1)).
\end{equation*}
The result follows by noticing that the last two maps $p \smas 1$ and $1 \smas \mu$ commute with each other, together with the fact that $s$ is a section of $p \smas 1$. For the third map, we can either do a similar diagram chase, or use the fact that $\Ct \smas \Ct$ is an $E_{\infty}$ ring spectrum, and so the third map is homotopic to the second map we just computed. The last map is forced to be nullhomotopic since
\begin{equation*}
\SCt \smas \SCt \simeq \SSSCt \vee \SSCt
\end{equation*}
and there are no non-trivial maps to both $\Ct$ and $\SCt$ by Corollary \ref{cor:mappingspacecoftau}.
\end{proof}
\end{lemma}

\vspace{0.3cm}

The additive splitting $\Ct \smas \Ct \simeq \Ct \vee \SCt$ gives the isomorphism
\begin{equation*}
\piss( \Ct \smas \Ct ) \cong \piss(\Ct) \oplus \btau \cdot \piss(\Ct).
\end{equation*}
The class $\btau$ has degree $|\btau| = (1,-1)$, and is the unit element of the shifted copy given by the composite
\begin{equation*}
S^{1,-1} \simeq S^{1,-1} \smas \Ss \stackrel{1 \smas i}{\ltoo} S^{1,-1} \smas \Ct \stackrel{s}{\lto} \Ct \smas \Ct.
\end{equation*}
We call it $\btau$ because it induces a $\tau$-Bockstein operations in $\HF \smas \Ct$-(co)homology, as we show in Propositions \ref{prop:HCtdualSteenrod} and \ref{prop:HCtSteenrod}. Lemma \ref{lem:multCtsmasCt} gives the following multiplicative description of the homotopy groups $\pi_{\ast,\ast}( \Ct \smas \Ct)$.

\begin{cor} \label{cor:multCtCt}
The $E_{\infty}$ ring spectrum $\Ct \smas \Ct$ has homotopy ring
\begin{equation*}
\pi_{\ast,\ast} \left( \Ct \smas \Ct \right) \cong \quotient{\piss \left( \Ct \right)[\btau]}{\btau^2},
\end{equation*}
where $|\btau| = (1,-1)$.
\end{cor}


\subsection{The Endomorphism Spectrum $\ECt$} \label{sec:EndCt}

In this Section we explicitly describe the homotopy type of $\ECt$ as a ring spectrum and give a presentation of its homotopy ring $\piss (\ECt)$, in the same way that we did for $\Ct \smas \Ct$. However, the endomorphism spectrum $\ECt$ is a little harder to understand than $\Ct \smas \Ct$. First, it is only an associative $A_{\infty}$ spectrum, whereas $\Ct \smas \Ct$ is $E_{\infty}$. Second, its multiplication comes from composition of morphisms and has nothing to do with the fact that $\Ct$ is a ring object, whereas the multiplication on $\Ct \smas \Ct$ is easy to describe in terms of the multiplication of $\Ct$. Finally, it turns out that out of the eight maps that assemble together to give the multiplication on $\ECt$, only three are forced to be nullhomotopic for degree reasons, whereas five where forced to be nullhomotopic for $\Ct \smas \Ct$. 

An important tool that we use is Spanier-Whitehead duality, adapted to the motivic setting from the categorical treatment in \cite[Chapter 3]{LMS}. We briefly recall some notation and elementary results from both \cite[Chapter 3]{LMS} and \cite[Sections 4.6-7]{HA}. Consider two motivic spectra $X$ and $Y$. If $X$ is dualizable, its Spanier-Whitehead dual is defined to be the motivic spectrum 
\begin{equation*}
DX \coloneqq F(X,\Ss).
\end{equation*}
In particular, finite cell complexes are dualizable. For spheres, there is a  canonical identification
\begin{equation} \label{eq:SWdualofS}
DS^{m,n} = F(S^{m,n}, \Ss) \simeq F(\Ss, S^{-m,-n}) \simeq S^{-m,-n}.
\end{equation}
Given a map $f \colon X \lto Y$ between dualizable motivic spectra, denote its Spanier-Whitehead dual by
\begin{equation*}
Df \coloneqq F(f,\Ss) \colon DY \lto DX.
\end{equation*}
If $X$ is dualizable, the smashing morphism $F(X, \Ss) \smas X \stackrel{\smas}{\lto} F(X,\Ss \smas X)$ is an equivalence, giving the equivalence
\begin{equation} \label{eq:equivSWdual}
DX \smas X = F(X, \Ss) \smas X \stackrel{\simeq}{\lto} F(X,\Ss \smas X) = \End(X).
\end{equation}
Denote the evaluation map that is adjoint to the identity map on $F(X,\Ss)$ by
\begin{equation*}
DX \smas X = F(X,\Ss) \smas X \stackrel{\ev}{\lto} \Ss.
\end{equation*}
The endomorphism spectrum $\End(X)$ is always a motivic $A_{\infty}$ ring spectrum with multiplication map given by the composite $\mu_{\End(X)}$ in the diagram
\begin{equation} \label{diag:multonEndX}
\begin{tikzpicture}
\matrix (m) [matrix of math nodes, row sep=2em, column sep=4em]
{ \End(X) \smas \End(X) &  & \End(X) \\
  DX \smas X \smas DX \smas X  &  DX \smas DX \smas X \smas X & DX \smas \Ss \smas X.  \\ };
\path[thick, -stealth]
(m-2-1) edge node[above] {$ 1 \smas \chi \smas 1 $} (m-2-2)
(m-2-2) edge node[above] {$ 1 \smas \ev \smas 1 $} (m-2-3);
\path[thick]
(m-1-1) edge[double, double distance=1.5pt] node[left] {can.} (m-2-1)
(m-1-3) edge[double, double distance=1.5pt] node[left] {can.} (m-2-3);
\path[thick, -stealth, dashed]
(m-1-1) edge node[above] {$ \mu_{\End(X)} $} (m-1-3);
\end{tikzpicture}
\end{equation}

The spectrum $\Ct$ is dualizable since it is a 2-cell complex. The $A_{\infty}$ ring structure on $\ECt$ can thus be understood in terms of Spanier-Whitehead duality. For this, we have to compute the homotopy type of the Spanier-Whitehead dual $D\Ct$ and identify the evaluation map $D\Ct \smas \Ct \stackrel{\ev}{\lto} \Ss$.

\begin{prop} \label{prop:SWdualofCt}
We have the following identifications.
\begin{enumerate}
\item The Spanier-Whitehead dual of $S^{0,-1} \stackrel{\tau}{\lto} \Ss$ is $D\tau \simeq \tau \colon \Ss \lto S^{0,1}$.
\item The Spanier-Whitehead dual of the cofiber sequence $$S^{0,-1} \stackrel{\tau}{\lto} \Ss \stackrel{i}{\lto} \Ct \stackrel{p}{\lto} S^{1,-1}$$ is the cofiber sequence
\begin{equation*}
S^{0,1} \stackrel{\tau}{\ltoback} \Ss \stackrel{p}{\ltoback} \CtS \stackrel{i}{\ltoback} S^{-1,1}.
\end{equation*}
In particular we have $Di \simeq p$ and $Dp \simeq i$, and a canonical (up to homotopy) identification 
\begin{equation} \label{eq:identifDCt}
D\Ct \simeq \CtS.
\end{equation}
\end{enumerate}
\begin{proof}
\mbox{}
\begin{enumerate}
\item Start with the map $S^{0,-1} \stackrel{\tau}{\lto} \Ss$. The functor $D=F(-, \Ss)$ and the canonical identification of equation \eqref{eq:SWdualofS} gives a map $\Ss \stackrel{D\tau}{\lto} S^{0,1}$, which by definition, sends $1$ to $\tau$ on $\pi_{0,0}$. Since it lives in the group $\left[ \Ss, S^{0,1} \right] \cong \hat{\Z}_2$ generated by $\tau$, we get that $D\tau \simeq \tau$.
\item Since the dualization functor $D$ preserves cofiber sequences, we get the cofiber sequence
\begin{equation*}
DS^{0,-1} \stackrel{D\tau}{\ltoback} D\Ss \stackrel{Di}{\ltoback} D\Ct \stackrel{Dp}{\ltoback} DS^{1,-1}.
\end{equation*}
To understand it, we use the canonical equivalences of equation \eqref{eq:SWdualofS} and embed it in the diagram
\begin{center}
\begin{tikzpicture}
\matrix (m) [matrix of math nodes, row sep=2em, column sep=4em]
{ DS^{0,-1} & D\Ss & D\Ct & DS^{1,-1}  \\
  S^{0,1}  &  \Ss & \CtS  & S^{-1,1}. \\ };
\path[thick, -stealth]
(m-1-2) edge node[above] {$ D\tau $} (m-1-1)
(m-1-3) edge node[above] {$ Di $} (m-1-2)
(m-1-4) edge node[above] {$ Dp $} (m-1-3)
(m-2-2) edge node[above] {$  \tau  $} (m-2-1)
(m-2-3) edge node[above] {$  p $} (m-2-2)
(m-2-4) edge node[above] {$  i $} (m-2-3);
\path[thick]
(m-1-1) edge[double, double distance=1.5pt] node[left] {can.} (m-2-1)
(m-1-2) edge[double, double distance=1.5pt] node[left] {can.} (m-2-2)
(m-1-4) edge[double, double distance=1.5pt] node[left] {can.} (m-2-4);
\path[thick, -stealth, dashed]
(m-2-3) edge (m-1-3);
\end{tikzpicture}
\end{center}
By the 5-lemma, the map $\CtS \lto D\Ct$ is an equivalence. Moreover, given two such equivalences, their difference would factor trough the map $p$ and thus trough $\Ss$. It follows that this equivalence is canonical up to homotopy, since by Corollary \ref{cor:htpycoftau} we have
\begin{equation*}
\pi_{0,0} (D\Ct) \cong \pi_{0,0} (\CtS) \cong \pi_{1,-1} (\Ct )= 0.
\end{equation*}
\end{enumerate}
\end{proof}
\end{prop}

\begin{lemma} \label{lem:evmaponCt}
Up to sign, the evaluation map $D\Ct \smas \Ct \stackrel{\ev}{\lto} \Ss$ is given by the commutative diagram
\begin{center}
\begin{tikzpicture}
\matrix (m) [matrix of math nodes, row sep=3em, column sep=4em]
{  D\Ct \smas \Ct & \Ss \\
   \CtS \smas \Ct & \CtS. \\};
\path[thick, -stealth]
(m-1-1) edge node[above] {ev} (m-1-2)
(m-2-1) edge node[above] {$ \mu $} (m-2-2)
(m-2-2) edge node[right] {$ p $} (m-1-2);
\path[thick]
(m-1-1) edge node[left] {$ \simeq $} node[right] {\emph{can.}} (m-2-1);
\end{tikzpicture}
\end{center}
\begin{proof}
We compute the abelian group of homotopy classes of maps $\left[ D\Ct \smas \Ct, \Ss \right]$. We have 
\begin{center}
\begin{tabular}{rlll}
$\left[ D\Ct \smas \Ct, \Ss \right]$   &$\cong  \left[ \CtS \smas \Ct, \Ss \right]$ & $\qquad$ & by equation \eqref{eq:identifDCt} \\
				  &$\cong  \left[ \CtS \vee \Ct, \Ss \right]$ & & by Lemma \ref{lem:equivCtsmasCt} \\
				  &$\cong  \left[ \CtS, \Ss \right] \oplus \left[ \Ct, \Ss \right]$ & & \\
				  &$\cong  \left[ \Ss, \Ss \right] \oplus 0$ & & via  $\CtS \stackrel{p}{\lto} \Ss$ \\
				  &$\cong \hat{\Z}_2$ & & 
\end{tabular}
\end{center}
which is generated by the identity. This means that $\left[ D\Ct \smas \Ct, \Ss \right]$ is generated by the composite
\begin{equation*}
D\Ct \smas \Ct \simeq \CtS \smas \Ct \stackrel{\mu}{\ltoo} \CtS \stackrel{p}{\lto} \Ss.
\end{equation*}
On the other side, by adjunction we have an isomorphism 
\begin{equation*}
\left[ D\Ct, D\Ct \right] \cong \left[ D\Ct \smas \Ct, \Ss \right],
\end{equation*}
which sends the identity map to the evaluation map (by definition of the evaluation map). This shows that $\ev$ is also one of the two units $\pm 1 \in \hat{\Z}_2$, finishing the proof. 
\end{proof}
\end{lemma}

\begin{lemma} \label{lem:multonECt}
Under the vertical identifications given by
\begin{center}
\begin{tikzpicture}
\matrix (m) [matrix of math nodes, row sep=2em, column sep=4em]
{  \ECt \smas \ECt  & \ECt  \\
   (\CtS \vee \Ct) \smas (\CtS \vee \Ct) & \CtS \vee \Ct \\
    (\CtS \smas \CtS) \vee (\CtS \smas \Ct) \vee (\Ct \smas \CtS) \vee (\Ct \smas \Ct) & \CtS \vee \Ct,  \\};
\path[thick, -stealth]
(m-1-1) edge node[above] {$ \mu_{\ECt} $} (m-1-2)
(m-2-1) edge (m-2-2);
\path[thick, -stealth, dashed]
(m-3-1) edge (m-3-2);
\path[thick]
(m-1-1) edge node[left] {$ \simeq $} (m-2-1)
(m-2-2) edge node[right] {$ \simeq $} (m-1-2)
(m-2-1) edge node[left] {$ = $} (m-3-1)
(m-2-2) edge node[right] {$ = $} (m-3-2);
\end{tikzpicture}
\end{center}
the multiplication on $\ECt$ is given by the maps
\begin{align*}
\CtS \smas \CtS &\stackrel{(p\smas 1,0)}{\ltooo} \CtS \vee \Ct \\
\CtS \smas \Ct &\stackrel{(\mu,0)}{\ltooo} \CtS \vee \Ct \\
\Ct \smas \CtS &\stackrel{(\mu, p\smas 1)}{\ltooo} \CtS \vee \Ct\\
\Ct \smas \Ct &\stackrel{(0,\mu)}{\ltooo} \CtS \vee \Ct.
\end{align*}
\begin{proof}[Sketch of proof]
This proof is by tedious diagram chases, and is in the spirit as the proof of Lemma \ref{lem:multCtsmasCt}. We will now briefly sketch the steps in the proof. The first part is to break $\ECt \smas \ECt$ in more manageable summands via Spanier-Whitehead duality, and the necessary identifications are done in Proposition \ref{prop:SWdualofCt}. We then use the definition of the multiplication map on $\ECt$ from diagram \eqref{diag:multonEndX}, as a composite of the factor swap map and the evaluation map. The evaluation map was explicitly computed in Lemma \ref{lem:evmaponCt}. The remainder of the proof consists on carefully identifying composites.
\end{proof}
\end{lemma}

The additive splitting $\ECt \simeq \Ct \vee \CtS$ gives the isomorphism
\begin{equation*}
\piss( \ECt ) \cong \piss(\Ct) \oplus \btau \cdot \piss(\Ct).
\end{equation*}
The class $\btau$ has degree $|\btau| = (-1,1)$, and is the unit element of the shifted copy given by the composite given by the composite
\begin{equation*}
\Ct \stackrel{p}{\lto} S^{1,-1} \stackrel{\Sigma i}{\lto} \Sigma^{1,-1} \Ct.
\end{equation*}
Lemma \ref{lem:multonECt} gives the following multiplicative description of the homotopy groups $\pi_{\ast,\ast}( \ECt)$.

\begin{cor} \label{cor:multEndCt}
The $A_{\infty}$ ring spectrum $\ECt$ has homotopy ring
\begin{equation*}
\piss \left( \ECt \right) \cong \quotient{\piss \left( \Ct \right)\langle \btau \rangle}{
\begin{array}{lr}
    \alpha \btau - (-1)^{|\alpha|}\btau \alpha = i\comp p(\alpha) \\
   \btau^2 = 0
  \end{array} 
  }
\end{equation*}
where $\btau$ is a non-commutative variable and $\alpha$ span the elements of $\piss(\Ct)$.
\end{cor}

\begin{remark}
The canonical inclusion $\Ct \lto \ECt$ is a map of $A_{\infty}$ ring spectra and on homotopy is the inclusion of $\piss (\Ct)$ onto the non-shifted factor. We can also think of the ring $\piss( \ECt)$ as being the abelian group 
\begin{equation*}
\piss( \ECt ) \cong \piss(\Ct) \oplus \btau \cdot \piss(\Ct)
\end{equation*}
with ring structure given by the following multiplication table
\begin{align*}
\alpha \comp \alpha' &= \alpha \alpha' \\
\alpha \comp \btau \alpha' &= (-1)^{|\alpha|} \btau \alpha \alpha' + (i \comp p(\alpha)) \alpha' \\
\btau \alpha \comp \alpha' &= \btau \alpha \alpha' \\
\btau \alpha \comp \btau \alpha' &= \btau (i\comp p(\alpha)) \alpha',
\end{align*}
where $\alpha, \alpha' \in \piss(\Ct)$ and $\btau\alpha, \btau\alpha' \in \btau \cdot \piss( \CtS)$.
\end{remark}

\begin{remark}
Since $\Ss \stackrel{i}{\lto} \Ct$ is the ring map which induces the $\pi_{\ast,\ast} (\Ss)$-module structure on $\pi_{\ast,\ast}( \Ct)$, we have the compatibility formula
\begin{equation*}
i(\alpha) \alpha' = \alpha \alpha' \qquad \qquad \text{for } \alpha \in \pi_{\ast,\ast} (\Ss), \alpha' \in \pi_{\ast,\ast} (\Ct).
\end{equation*}
The first multiplication uses the ring structure of $\Ct$ while the second uses the $\Ss$-module structure on $\Ct$. This simplifies some of the formulas of Corollary \ref{cor:multEndCt}, for example by $\btau \alpha \comp \btau \alpha' = \btau p(\alpha) \alpha'$ since $p(\alpha)$ is in the homotopy groups of the motivic sphere. 
\end{remark}

\section{Examples of $\Ct$-modules} \label{sec:Ctmodules}

Since the 2-cell complex $\Ct$ is a (cofibrant) commutative ring spectrum, we can use \cite[Section 2.8]{Pelaez} to endow the category $\CtMod$ with a closed symmetric monoidal model structure. The closed monoidal structure is given by the the relative smash product $- \smasCt -$ and the internal function spectrum $\FCt(-,-)$. Moreover, the model structure is created by the forgetful functor, and is thus part of the Quillen adjunction
\begin{equation} \label{eq:freeforgetCt}
\SptC = {}_{\Ss}\textbf{Mod} \overunder{- \smas \Ct}{U}{\adj} \CtMod.
\end{equation}
In this Section we will first give some elementary lemmas about the category $\CtMod$, and then study some important spectra that are induced up from $\Ss$-modules by smashing with $- \smas \Ct$. We call such a spectrum a \emph{$\Ct$-induced spectrum}.

We start with the $\Ct$-induced Eilenberg-Maclane spectrum $\HF \smas \Ct$ which has homotopy groups $\piss( \HF \smas \Ct) \cong \F_2$ in degree $(0,0)$. We will compute its Steenrod algebra of operations (and its dual) as a Hopf algebra, both in $\SptC$ and $\CtMod$. This computation is used in future work \cite{GheKwn} to construct Morava $K$-theories for the motivic $w_i$ periodic operators. The first operator $w_1$ was introduced in \cite{Andrews}. We then show that the $\Ct$-induced Moore spectrum $\Stwotau$ admits a unique structure of an $E_{\infty}$ algebra over $\Ct$. We also observe that it admits a $v_1^1$-self map, whereas $\Ss/2$ only admits a $v_1^4$-self map. Finally, we compute the homology and homotopy of the $\Ct$-induced connective algebraic and hermitian $K$-theory spectra $\kgl$ and $\kq$. Here again an interesting phenomenon arises in hermitian $K$-theory: an obstruction is killed and we can see the element $v_1^2$ in the homotopy of $\kq \smas \Ct$, whereas we only see its square $v_1^4$ in $\kq$.

\subsection{Elementary Results on $\Ct$-Modules} \label{sec:GenCtmod}

Let $X$ be a (left) $\Ct$-module with action map $\phi_X \colon \Ct \smas X \lto X$. The left unitality condition says that the triangle in the diagram
\begin{center}
\begin{tikzpicture}
\matrix (m) [matrix of math nodes, row sep=3em, column sep=3em]
{  S^{0,-1} \smas X & \Ss \smas X & \Ct \smas X & S^{1,-1} \smas X \\
                    &             &           X &  \\};
\path[thick, -stealth, font=\small]
(m-1-1) edge node[above] {$ \tau $} (m-1-2)
(m-1-2) edge node[auto] {$ i $} (m-1-3)
(m-1-3) edge node[auto] {$ p $} (m-1-4)
(m-1-2) edge node[auto] {$  \simeq $} (m-2-3)
(m-1-3) edge node[right] {$  \phi_X $} (m-2-3);
\end{tikzpicture}
\end{center}
commutes, i.e., that $\phi_X$ is a retraction of the unit. This produces a splitting 
\begin{equation} \label{eq:CtsmasXsplits}
\Ct \smas X \stackrel{(\phi_X,p)}{\ltooo} X \vee \Sigma^{1,-1} X
\end{equation}
up to homotopy, whose inverse map requires a choice of section of $p$. There is however a canonical choice of section given by the composite
\begin{equation*}
S^{1,-1} \smas X = S^{1,-1} \smas \Ss \smas X \stackrel{\id \smas i \smas \id}{\ltooo} S^{1,-1} \smas \Ct \smas X \stackrel{s \smas \id}{\lto} \Ct \smas \Ct \smas X \stackrel{\id \smas \phi_X}{\ltoo} \Ct \smas X,
\end{equation*}
by using the canonical section $s \colon \SCt \lto \Ct \smas \Ct$ from Lemma \ref{lem:equivCtsmasCt}. The Betti realization functor $\SptC \lto \Spt$ naturally extends to $\CtMod$ by composing with the forget functor
\begin{equation*}
\CtMod \lto \SptC \stackrel{\Betti}{\lto} \Spt.
\end{equation*}


\begin{lemma} \label{lem:BettiRelofCtMod}
Every $\Ct$-module realizes to a contractible spectrum in $\Top$.
\begin{proof}
Consider a spectrum $X \in \SptC$ endowed with a structure of $\Ct$-module. Since the Betti realization functor is (strict) symmetric monoidal and sends $\Ct$ to a contractible spectrum, we have
\begin{equation*}
\Betti( \Ct \smas X) \simeq \Betti(\Ct) \smas \Betti(X) \simeq \ast.
\end{equation*}
It follows that $\Betti(X) \simeq \ast$ as $X$ is a retract of $\Ct \smas X$ by equation \eqref{eq:CtsmasXsplits}.
\end{proof}
\end{lemma}

The next two elementary lemmas will often be used for studying $\Ct$-induced spectra.

\begin{lemma} \label{lem:htpyhlgyoftaufree}
Let $X$ be a spectrum with $\tau$-free homotopy (resp. homology) groups, i.e., multiplication by $\tau$ is injective on $\piss(X)$ (resp. on $\HFho(X)$). Then the homotopy (resp. homology) groups of the $\Ct$-induced spectrum $X \smas \Ct$ are given by
\begin{equation*}
\piss(X \smas \Ct) \cong \quotient{\piss(X)}{\tau}  \qquad \text{ (resp. } \HFho(X \smas \Ct) \cong \quotient{\HFho(X)}{\tau} \text{ ).}
\end{equation*}
Moreover if $X$ is an $E_{\infty}$ ring spectrum, then this isomorphism is a ring isomorphism.
\begin{proof}
This follows by the long exact sequence induced from the cofiber sequence
\begin{equation*}
\Sigma^{0,-1} X \stackrel{\tau}{\lto} X \stackrel{i}{\lto} \Ct \smas X
\end{equation*}
since multiplication by $\tau$ is injective. Moreover, if $X$ is an $E_{\infty}$ ring spectrum, then the map 
\begin{equation*}
S^{0,0} \smas X \stackrel{i \smas \id}{\ltoo} \Ct \smas X
\end{equation*}
is a map of $E_{\infty}$ ring spectra as well.
\end{proof}
\end{lemma}

\begin{lemma} \label{lemma:cohomoftaufree}
Let $X$ be a spectrum with $\tau$-free $\HF$-cohomology groups, i.e., multiplication by $\tau$ is injective on $\HFco(X)$. Then the cohomology groups of the $\Ct$-induced spectrum $X \smas \Ct$ are given by
\begin{equation*}
\HFco(X \smas \Ct) \cong \quotient{\HFco( \Sigma^{1,-1}X)}{\tau}.
\end{equation*}
\begin{proof}
Similarly to the proof of Lemma \ref{lem:htpyhlgyoftaufree}, this just follows by the long exact sequence induced from the cofiber sequence
\begin{equation*}
\Ct \smas X  \lto \Sigma^{1,-1}X \stackrel{\tau}{\lto} \Sigma^{1,0}X
\end{equation*}
since multiplication by $\tau$ is injective. 
\end{proof}
\end{lemma}

\subsection{The $\Ct$-Induced Eilenberg-Maclane Spectrum} \label{sec:HFtwomodtau}

Consider the $\Ct$-induced Eilenberg-Maclane spectrum 
\begin{equation*}
\HCt \coloneqq \HF \smas \Ct,
\end{equation*}
which has homotopy $\piss (\HCt) \cong \F_2$ concentrated in degree $(0,0)$ by Lemma \ref{lem:htpyhlgyoftaufree}. Unlike $\HF$, this spectrum detects both cells of $\Ct$ since
\begin{equation*}
\HFco( \Ct ) \cong 
  \begin{cases} 
   \F_2 & \text{if } (\ast,\ast) = (0,0) \\
   \F_2 & \text{if } (\ast,\ast) = (1,-1) \\
   0       & \text{otherwise.}
  \end{cases}
\end{equation*}
This spectrum plays an important role in the theory of motivic periodicities, as it is the building block of the Morava K-theories $K(w_n)$ and of the Brown-Peterson spectrum $wBP$ that we construct in \cite{GheKwn}. Denote the $\HCt$-Steenrod algebra of operations in $\HCt$-cohomology by
\begin{equation*}
\ACtau \cong \pi_{-\ast,-\ast} \left( F(\HCt, \HCt) \right),
\end{equation*}
and its dual algebra of co-operations in $\HCt$-homology by
\begin{equation*}
\ACtaudual \cong \pi_{\ast,\ast} \left( \HCt \smas \HCt \right) .
\end{equation*} 
The two main ingredients for these computations are our previous knowledge of the $\HF$-Steenrod algebra $\AC$, which we recalled in Section \ref{sec:SteenrodAlg}, and the descriptions of $\Ct \smas \Ct$ and $\ECt$ from Section \ref{sec:computationsCt}. Since $\tau \in \M$ is an element of the base ring, there is an induced Hopf algebra structure over $\M / \tau \cong \F_2$ on the quotients $\quotient{\AC}{\tau}$ and $\quotient{\ACdual}{\tau}$. 

\begin{prop} \label{prop:HCtdualSteenrod}
The dual $\HCt$-Steenrod algebra $\ACtaudual$ has the following Hopf algebra structure
\begin{equation*}
\ACtaudual \cong \quotient{\ACdual}{\tau} \otimes E(\btau) \cong \F_2[\xi_1, \xi_2, \ldots ] \otimes E(\tau_0, \tau_1, \ldots) \otimes E(\btau)
\end{equation*}
where $\btau$ is a $\tau$-Bockstein in degree $(1,-1)$ which is primitive in the coalgebra structure.
\begin{proof}
The dual $\HCt$-Steenrod algebra is given by the homotopy groups of the $E_{\infty}$ ring spectrum
\begin{equation*}
\HCt \smas \HCt = \HF \smas \Ct \smas \HF \smas \Ct \simeq \HF \smas \HF \smas \Ct \smas \Ct.
\end{equation*}
Since $\piss(\HCt) \cong \F_2$, the left and right units of the Hopf algebroid $\piss\left( \HCt \smas \HCt \right)$ are flat maps and they agree, turning it into a Hopf algebra. If we smash the canonical equivalence $\Ct \smas \Ct \simeq \Ct \vee \SCt$ of Lemma \ref{lem:equivCtsmasCt} with $\HF \smas \HF$, we get an additive splitting
\begin{equation*}
\HCt \smas \HCt \simeq \left( \HF \smas \HF \smas \Ct \right)  \vee \left( \Sigma^{1,-1} \HF \smas \HF \smas \Ct \right),
\end{equation*}
into two wedge summands that we can understand individually. Since the dual Steenrod algebra $\ACdual$ is $\tau$-free, Lemma \ref{lem:htpyhlgyoftaufree} gives a ring description of the homotopy
\begin{equation*}
\piss( \HF \smas \HF \smas \Ct) \cong \quotient{\ACdual}{\tau},
\end{equation*}
and thus the dual $\HCt$-Steenrod algebra is a free module of rank 2 over $\ACdual$. The first generator in degree $(0,0)$ is the unit given by the ring map
\begin{equation*}
\Ss \stackrel{i}{\ltoo} \HCt \smas \HCt.
\end{equation*}
The second generator in degree $(1,-1)$ that we call $\btau$ is given by the map 
\begin{equation*}
\btau \colon S^{1,-1} \stackrel{ i}{\lto} \Sigma^{1,-1} \Ct \stackrel{s}{\lto} \Ct \smas \Ct \stackrel{i \smas i}{\lto} \HCt \smas \HCt,
\end{equation*}
where $i$ denotes the inclusion of the bottom cell and $s$ denotes the canonical section of $\mu$, as in Lemma \ref{lem:equivCtsmasCt}. We choose the name $\btau$ because its dual element in the $\HCt$-Steenrod algebra does behave like a $\tau$-Bockstein in cohomology, as  we explain in Proposition \ref{prop:HCtSteenrod}. To finish the description of the ring structure of $\ACtaudual$, we have to compute the product $\btau \cdot \btau$ which lands in degree $(2,-2)$. This product is the homotopy class of the composite
\begin{equation*}
\btau \cdot \btau \colon S^{1,-1} \smas S^{1,-1} \stackrel{\btau \smas \btau}{\ltooo} \HCt \smas \HCt \smas \HCt \smas \HCt \stackrel{\mu}{\lto} \HCt \smas \HCt
\end{equation*}
which is nullhomotopic since $\mu_{\Ct} \comp s \simeq 0$. This gives the ring structure as the tensor products
\begin{equation*}
\ACtaudual \cong \ACdual \otimes E(\btau) \cong \F_2[\xi_1, \xi_2, \ldots] \otimes E(\tau_0, \tau_1, \ldots) \otimes E(\btau).
\end{equation*}
For the coalgebra structure, the counit is forced as there is only a copy of $\F_2$ in degree $(0,0)$. It thus only remains to compute the coproduct. The ring map
\begin{equation*}
\HF \stackrel{i}{\lto} \HCt
\end{equation*}
induces the following map of Hopf algebras 
\begin{equation*}
\ACdual \stackrel{\psi}{\lto} \ACtaudual \cong \ACdual / \tau \otimes E(\btau) \colon a \lmapsto a \otimes 1,
\end{equation*}
which can be factored as reduction modulo $\tau$ and then inclusion into the $- \otimes 1$ factor. It follows that the coproduct $\Delta(a \otimes 1)$ can be computed by choosing a pre-image $a$ of $a \otimes 1$, computing the coproduct in $\ACdual$, and then pushing it back via $\psi$. Since the coproduct formula on the $\xi_i$'s and $\tau_i$'s in $\ACdual$ does not involve any $\tau$-multiples, the exact same formula holds for the coproduct of elements of the form $a \otimes 1 \in \ACtaudual$. It only remains to compute the diagonal on the element $1 \otimes \btau$. We show in the next Proposition \ref{prop:HCtSteenrod} that its dual is exterior in the algebra structure of $\ACtau$, implying that $1 \otimes \btau$ is primitive.
\end{proof}
\end{prop}

\begin{prop} \label{prop:HCtSteenrod}
The $\HCt$-Steenrod algebra $\ACtau$ has the following Hopf algebra structure
\begin{equation*}
\ACtau \cong \quotient{\AC}{\tau} \otimes E(\btau) 
\end{equation*}
where $\btau$ is a $\tau$-Bockstein in degree $(1,-1)$ which is primitive in the coalgebra structure.
\begin{proof}
Since $\Ct$ is dualizable we can rewrite
\begin{equation*}
F( \HCt, \HCt) = F( \HF \smas \Ct, \HF \smas \Ct) \simeq F(\HF, \HF) \smas \Ct \smas \DCt.
\end{equation*}
By the identification of Section \ref{sec:EndCt} we further have
\begin{equation*}
F( \HCt, \HCt) \simeq \left( F(\HF,\HF) \smas \Ct \right) \vee \left( \Sigma^{-1,1} F(\HF, \HF) \smas \Ct \right).
\end{equation*}
By Lemma \ref{lem:htpyhlgyoftaufree} we get that $\ACtau$ is a free $\AC/\tau$-module of rank 2 with generators given by the operations
\begin{equation*}
\id \colon \HCt \lto \HCt \qquad \text{ and } \qquad \btau \colon \HCt \stackrel{p}{\lto} \Sigma^{1,-1} \HF \stackrel{i}{\lto} \Sigma^{1,-1} \HCt,
\end{equation*}
where $p$ denotes the projection of $\Ct$ on its top cell, while $i$ denotes the inclusion of it bottom cell. The definition of $\btau$ explains why we call it a $\tau$-Bockstein. Since the Steenrod algebra is defined as negative homotopy groups of the endomorphism spectrum, the $\tau$-Bockstein $\btau$ is in degree $(1,-1)$. This settles the additive structure of $\ACtau$, and it remains to understand its Hopf algebra structure. Since $\ACtaudual$ is a Hopf algebra of finite type, we can dualize its structure from Proposition \ref{prop:HCtdualSteenrod} to get the desired Hopf algebra structure of $\ACtau$. Recall that we did not yet finish the proof of Proposition \ref{prop:HCtdualSteenrod}, as we still have to show that $\btau \in \ACtaudual$ is primitive. This is equivalent to $\btau \in \ACtau$ being exterior, which is clear since it is the composite
\begin{equation*}
\btau \comp \btau \colon \HCt \stackrel{p}{\lto} \Sigma^{1,-1} \HF \stackrel{i}{\lto} \Sigma^{1,-1} \HCt \stackrel{p}{\lto} \Sigma^{1,-1} \HF \stackrel{i}{\lto} \Sigma^{1,-1} \HCt,
\end{equation*}
which is nullhomotopic as $p \comp i \simeq 0$.
\end{proof}
\end{prop}

\begin{remark}[$\Ct$-linear $\HCt$-homology and cohomology]
We can define the $\Ct$-linear homology and cohomology of a $\Ct$-module $X$ to be
\begin{equation*}
\HCt_{\ast,\ast}^{\Ct}(X) \coloneqq \piss(\HCt \smasCt X) \qquad \text{ and } \qquad \HCt^{\ast,\ast}_{\Ct}(X) \coloneqq \pi_{-\ast,-\ast} \left( \FCt(X,\HCt) \right).
\end{equation*}
The relevant $\HCt$-Steenrod algebra of $\Ct$-linear operations and co-operations are then
\begin{equation*}
\pi_{-\ast,-\ast} \left( \FCt(\HCt, \HCt) \right) \qquad \text{ and } \qquad \ACtaudual \cong \pi_{\ast,\ast} \left( \HCt \smasCt \HCt \right).
\end{equation*} 
Their computation follows from Lemmas \ref{lem:htpyhlgyoftaufree} and \ref{lemma:cohomoftaufree}, and the result is the usual motivic Steenrod algebra and its dual, modulo $\tau$. The only difference with the computations of Propositions \ref{prop:HCtdualSteenrod} and \ref{prop:HCtSteenrod} is that the $\Ct$-linear Steenrod algebras do not contain the $\tau$-Bockstein element $\beta_{\tau}$. In particular, the dual $\Ct$-linear $\HCt$-Steenrod algebra enjoys the nice formula
\begin{equation*}
\F_2[\xi_1, \xi_2, \ldots ] \otimes E(\tau_0, \tau_1, \ldots) 
\end{equation*}
that is very reminiscent of the odd-primary classical Steenrod algebra.
\end{remark}


\subsection{The $\Ct$-Induced Moore Spectrum}

Denote by $S^0/2$ the mod 2 Moore spectrum in the usual category of topological spectra $\Spt$. Recall that the classical Toda bracket $\langle 2, \eta, 2 \rangle = \eta^2$ implies that $\pi_2(S^0/2) \cong \Z/4$. This shows that multiplication by 2 is not a nullhomotopic map on $S^0/2$, and thus that there is no possible filler in the diagram
\begin{center}
\begin{tikzpicture}
\matrix (m) [matrix of math nodes, row sep=3em, column sep=4em]
{ S^0 \smas S^0/2 & S^0 \smas S^0/2 & S^0/2 \smas S^0/2 & \Sigma^1 S^0/2 \\
  & & S^0/2. & \\};
\path[thick, -stealth, font=\small]
(m-1-1) edge node[above] {$ 2 $} (m-1-2)
(m-1-2) edge (m-1-3)
(m-1-3) edge (m-1-4)
(m-1-2) edge node[auto] {$ \simeq $} (m-2-3);
\path[thick, -stealth, dashed, font=\small]
(m-1-3) edge node[right] {$ \nexists \ \mu $} (m-2-3);
\end{tikzpicture}
\end{center}
This shows that there exists no left unital multiplication on $S^0/2$.

Denote now the motivic mod 2 Moore spectrum by $\Ss/2$. Similarly, we can compute the motivic homotopy group $\pi_{2,0} (\Ss/2) \cong \Z/4$ via the same argument. More precisely, the analoguous Toda bracket is $\langle 2, \tau \eta, 2 \rangle = \tau^2 \eta^2$, where $\eta \in \pi_{1,1} (\Ss)$ and thus $\tau \eta \in \pi_{1,0} (\Ss)$. This again implies that there is no left unital multiplication on the Moore spectrum $\Ss/2$. Observe that this could also have been noticed by the fact that a left unital multiplication on $\Ss/2$ would induce one on $S^0/2$ by Betti realization. 

Denote the cofiber of multiplication by $\tau$ on $\Ss/2$ by $\Stwotau$. This spectrum does admit a left unital multiplication since 
\begin{equation*}
\langle 2, \eta, 2 \rangle = \tau \eta^2 \equiv 0 \qquad \text{ modulo } \tau.
\end{equation*}
This does not imply that there is a ring structure on $\Stwotau$ as this bracket is just one possible obstruction (the obstruction to left unitality). In Theorem \ref{thm:Smod2tauisEoo} we show that all obstructions are of this type and that $\Stwotau$ admits the structure of an $E_{\infty}$ algebra over $\Ct$. 

Since cofibers in $\Ct$-modules can be computed in the underlying category of motivic spectra, it follows that the cofiber of 2 on $\Ct$ has underlying spectrum $\Stwotau$. Consider now $\Stwotau$ as a $\Ct$-module, for example as constructed in the category $\CtMod$ by the cofiber sequence 
\begin{equation} \label{eq:cofseqStwotau}
\Ct \stackrel{2}{\lto} \Ct \stackrel{i}{\lto} \Stwotau \stackrel{p}{\lto} \Sigma^{1,0}\Ct.
\end{equation}
To equip $\Stwotau$ with an $E_{\infty}$ $\Ct$-algebra structure, we will proceed very similarly as in Section \ref{sec:CtEinfty}, which we refer to for more details.

\begin{prop}  \label{prop:Stwotauishtpyunital}
There is a unique homotopy unital and homotopy commutative $\Ct$-algebra structure on $\Stwotau$.
\begin{proof}
The computation of $\left[ \Stwotau, \Stwotau \right]_{\Ct} \cong \Z/2$ generated by the identity map shows that $\cdot 2$ is nullhomotopic on $\Stwotau$, providing a left unital multiplication $\mu$ from the diagram
\begin{center}
\begin{tikzpicture}
\matrix (m) [matrix of math nodes, row sep=3em, column sep=4em]
{ \Ct \smasCt \Stwotau & \Ct \smasCt \Stwotau & \Stwotau \smasCt \Stwotau & \Sigma^{1,0} \Ct \smasCt \Stwotau \\
  & & \Stwotau. & \\};
\path[thick, -stealth, font=\small]
(m-1-1) edge node[above] {$ 2 $} (m-1-2)
(m-1-2) edge node[auto] {$ i_L $}  (m-1-3)
(m-1-3) edge node[auto] {$ p_L  $} (m-1-4)
(m-1-2) edge node[auto] {$ \simeq $} (m-2-3);
\path[thick, -stealth, dashed, font=\small]
(m-1-3) edge node[right] {$ \exists \ \mu $} (m-2-3);
\end{tikzpicture}
\end{center}
The computation $\left[ \Sigma^{1,0} \Ct \smasCt \Stwotau, \Stwotau \right]_{\Ct} = 0$ shows that there is a unique left unital multiplication up to homotopy on $\Stwotau$. As in Lemma \ref{lem:equivCtsmasCt}, it also implies that there is a unique section $s$ of $p_L$, giving a canonical additive splitting
\begin{equation} \label{eq:splittingStwotausmashitself}
\Stwotau \smasCt \Stwotau \simeq \Stwotau \vee \Sigma^{1,0} \Stwotau.
\end{equation}
The induced multiplication $\widetilde{\mu}$ after this identification is again just projection onto the first factor, and the factor swap map $\chi$ is given by the following diagram
\begin{center}
\begin{tikzpicture} [ampersand replacement=\&]
\matrix (m) [matrix of math nodes, row sep=3em, column sep=4em]
{  \Stwotau \smasCt \Stwotau \& \Stwotau \smasCt \Stwotau \\
   \Stwotau \vee \Sigma^{1,0} \Stwotau \& \Stwotau \vee \Sigma^{1,0} \Stwotau. \\};
\path[thick, -stealth]
(m-1-1) edge node[above] {$ \chi $} (m-1-2)
(m-2-1) edge node[left] {$ i_L + s $} (m-1-1)
(m-1-2) edge node[right] {$ ( \mu, p_L) $} (m-2-2);
\path[thick, -stealth, dashed]
(m-2-1) edge node[above] {$ \left[ \begin{smallmatrix} 1&0\\ i \comp p&1 \end{smallmatrix}  \right] $} (m-2-2);
\end{tikzpicture}
\end{center}
The matrix can be completely determined since $\left[ \Stwotau, \Stwotau \right]_{\Ct} \cong \Z/2$. By an easy matrix multiplication as in Proposition \ref{prop:Ctismonoidinhtpy}, this shows that $\mu$ is right unital and homotopy commutative.
\end{proof}
\end{prop}

The next step is to show that this (unique) multiplication map $\mu$ on $\Stwotau$ can be extended to an $E_{\infty}$ multiplication. We proceed in the exact same way as we did in Proposition \ref{prop:multhtpyassoc} and Theorem \ref{thm:multCthtpycom}.

\begin{thm} \label{thm:Smod2tauisEoo}
The $\Ct$-algebra structure on $\Stwotau$ can be uniquely extended to an $E_{\infty}$ structure.
\begin{proof}
We first extend it to an $A_{\infty}$ structure as in Proposition \ref{prop:multhtpyassoc}, with obstructions living in the abelian group
\begin{equation*}
\left[ \Sigma^{n-3,0} (\Sigma^{1,0}\Ct)^{\smas n}, \Stwotau \right]_{\Ct} \cong \left[ \Sigma^{2n-3,0} \Ct^{\smas n}, \Stwotau \right]_{\Ct}
\end{equation*}
for $n \geq 3$. Here we used $\Sigma^{1,0}\Ct$ since it is the cofiber of the unit map $\Ct \stackrel{i}{\lto} \Stwotau$. By using the decomposition formula for $\Ct^{\smas n}$ from Corollary \eqref{cor:decCtsmas}, the obstructions live in the group
\begin{equation*}
\bigoplus_{i=0}^{n} \binom{n}{i} \left[ \Sigma^{2n -3 + i,-i} \Ct, \Stwotau \right]_{\Ct}.
\end{equation*}
By the free-forget adjunction these groups are
\begin{equation*}
\pi_{2n -3 + i,-i} (\Stwotau).
\end{equation*}
For $n \geq 3$ and for any $0 \leq i \leq n$ this homotopy group is zero, making the obstruction group zero and allowing $\mu$ to extend to an $A_{\infty}$ structure. Similarly the obstructions for uniqueness live in zero groups, showing that $\Stwotau$ admits a unique $A_{\infty}$ algebra structure over $\Ct$.

The $A_3$ structure gives an associative homotopy, and thus we now have a unital, associative and commutative monoid in the homotopy category. This is a 3-stage in Robinsin's obstruction theory, so we can apply Corollary \ref{cor:EooObstTheoryMachinery} to extend it to an $E_{\infty}$ ring structure. The obstructions live in 
\begin{equation*}
\left[ \Sigma^{n-3,0} \Stwotau^{\smas m}, \Stwotau \right]_{\Ct}
\end{equation*}
for $n \geq 4$ and $2 \leq m \leq n$, where the smash product is over $\Ct$. As in the proof of Theorem \ref{thm:multCthtpycom}, we first break the source in smaller pieces by recursively using equation \eqref{eq:splittingStwotausmashitself}. It is then easy to show that all of those groups are zero by using cofiber sequences in the first variable to reduce it to homotopy groups of $\Stwotau$. Similarly, the obstructions for uniqueness live in
\begin{equation*}
\left[ \Sigma^{n-2,0} \Stwotau^{\smas m}, \Stwotau \right]_{\Ct}
\end{equation*}
for $n \geq 4$ and $2 \leq m \leq n$. We show by the exact same method that all those groups are zero, finishing the proof.
\end{proof}
\end{thm}

\begin{remark} \label{rem:exwhyalgcatarebetter}
The fact that multiplication by $2$ is nullhomotopic on $\Stwotau \simeq \Ct/2$ is not so surprising, as $\Ct$ is of somehow of algebraic nature. In fact, multiplication by $n$ on $X/n$ is always nullhomotopic in such algebraic categories, as explained in \cite[Proposition 1]{Schwede}.
\end{remark}

\begin{remark} \label{rem:v1onSmodtwotau}
The Toda bracket $\langle 2, \eta, 2 \rangle = \eta^2$ is also responsible for the non-existence of a $v_1^1$-self map on the topological Moore spectrum $S^0/2$. This is illustrated in the diagram
\begin{center}
\begin{tikzpicture}
\matrix (m) [matrix of math nodes, row sep=3em, column sep=4em]
{  S^2/2 & S^2 &  S^2  \\
   S^0/2 & S^1 &  S^1. \\};
\path[thick, -stealth, font=\small]
(m-1-2) edge node[above] {$ i $} (m-1-1)
(m-1-3) edge node[above] {$ 2 $}  (m-1-2)
(m-2-1) edge node[above] {$ p $} (m-2-2)
(m-2-2) edge node[above] {$ 2 $} (m-2-3)
(m-1-2) edge node[right] {$ \eta $} (m-2-2);
\path[dashed, -stealth, font=\small]
(m-1-2) edge node[above=5pt, left] {$ \exists \ \widetilde{\eta} $} (m-2-1)
(m-1-1) edge node[left] {$ \nexists $} (m-2-1);
\end{tikzpicture}
\end{center}
The map $\widetilde{\eta}$ exists since $2\eta = 0$, but there is no $v_1^1$-self map as $2 \cdot \widetilde{\eta} \neq 0$.  Motivically, the same diagram has the same problem because of the non-vanishing of the bracket $\langle 2, \eta, 2 \rangle = \tau \eta^2$. However, in $\Ct$-modules this bracket vanishes and the $\Ct$-induced Moore spectrum admits a $v_1^1$-self map. The diagram
\begin{center}
\begin{tikzpicture}
\matrix (m) [matrix of math nodes, row sep=3em, column sep=4em]
{  \Sigma^{2,1}\Stwotau &  \Sigma^{2,1}\Ct &   \Sigma^{2,1}\Ct \\
   \Stwotau & \Sigma^{1,0}\Ct &  \Sigma^{1,0}\Ct \\};
\path[thick, -stealth, font=\small]
(m-1-2) edge node[above] {$ i $} (m-1-1)
(m-1-3) edge node[above] {$ 2 $}  (m-1-2)
(m-2-1) edge node[above] {$ p $} (m-2-2)
(m-2-2) edge node[above] {$ 2 $} (m-2-3)
(m-1-2) edge node[right] {$ \eta $} (m-2-2);
\path[dashed, -stealth, font=\small]
(m-1-2) edge node[above=5pt, left] {$ \exists \ \widetilde{\eta} $} (m-2-1)
(m-1-1) edge node[left] {$ \exists \ v_1 $} (m-2-1);
\end{tikzpicture}
\end{center}
exhibits this $v_1$-self map
\begin{equation*}
\Sigma^{2,1} \Stwotau \stackrel{v_1}{\lto} \Stwotau.
\end{equation*}
More precisely, this follows since the computation $\left[ \Sigma^{2,1} \Ct, \Stwotau \right] \cong \Z/2$ forces the relation $2 \cdot \tilde{\eta} \simeq 0$.
\end{remark}


\subsection{The $\Ct$-Induced connective Algebraic and Hermitian $K$-Theory Spectra}

Consider the motivic algebraic $K$-theory spectrum $KGL$ constructed in \cite{Voe}. This spectrum represents algebraic $K$-theory on schemes. More precisely, given any scheme $X$, the $KGL$-cohomology of its stabilization $\Sigma^{\infty}_{+} X$ computes the algebraic $K$-theory of the scheme $X$. Consider now its connective cover $\kgl$ as described in \cite{IS} over $\Spec \C$ and in \cite{NSO} over more general basis. It is shown in \cite{NSO} that both $KGL$ and $kgl$ admit a unique $E_{\infty}$ ring structure. Recall that we work in the 2-completed category, and we use $\kgl$ to denote the 2-completed connective algebraic $K$-theory spectrum. Its coefficients and mod 2 homology of $\kgl$ over $\Spec \C$ are computed in \cite{IS} and given by
\begin{equation*}
\piss(\kgl) \cong \hat{\Z}_2[\tau, v_1] \qquad \text{ and } \qquad \HFho( \kgl) \cong \quotient{ \F_2[\tau][\xi_1,\xi_2, \ldots][\tau_2, \tau_3,\ldots]}{\tau_i^2 = \tau \xi_{i+1}},
\end{equation*}
where the element $v_1$ is in degree $(2,1)$ and corresponds to the usual Bott periodicity. Its homology is written as a subalgebra of the mod 2 homology of $\HF$ recalled in equation \eqref{eq:dualSteenrodalg}.

\begin{figure}[h] 
\begin{tikzpicture}[scale=0.5] 
\centering

\draw [very thin, -stealth] (0,-0.5) grid (0,8);
\node [left= 4pt] at (0,8) { filtration };
\draw [very thin, -stealth] (-0.5,0) -- (12,0);
\node [below=4pt] at (12,0) { stem };
\draw[step=1cm,gray,very thin, dotted] (0,0) grid (11.9,7.9);
\node [below=4pt] at (4,0) { $4$ };
\node [below=4pt] at (8,0) { $8$ };

\small{
\foreach \x in {0,...,7} { \draw[fill] (0,\x) circle [radius=0.1];  } 
\draw [-stealth] (0,0) -- (0,8);
						   			   
\foreach \x in {5,...,7} { \draw[fill] (4,\x) circle [radius=0.1]; } 
\draw[fill] (4,3) circle [radius=0.1];
\draw[fill] (4.2,4) circle [radius=0.1];
\draw (4,3) -- (4.2,4) -- (4,5);
\draw [-stealth] (4,5) -- (4,8);

\foreach \x in {4,...,7} { \draw[fill] (8,\x) circle [radius=0.1]; } 
\draw [-stealth] (8,4) -- (8,8);

\draw[fill] (1,1) circle [radius=0.1];
\draw[fill] (2,2) circle [radius=0.1];
\draw (0,0) -- (1,1)-- (2,2)-- (3,3);
\draw [-stealth, red] (3,3) -- (6,6);
\node[mark size=1.5pt,color=red, right] at (3-0.88,3) {\pgfuseplotmark{triangle*}};
\node[mark size=1.5pt,color=red, right] at (4-0.88,4) {\pgfuseplotmark{triangle*}};
\node[mark size=1.5pt,color=red, right] at (5-0.88,5) {\pgfuseplotmark{triangle*}};

\draw (8,4) -- (11,7);
\draw [-stealth, red] (11,7) -- (12,8);
\draw[fill] (9,5) circle [radius=0.1];
\draw[fill] (10,6) circle [radius=0.1];
\node[mark size=1.5pt,color=red, right] at (11-0.88,7) {\pgfuseplotmark{triangle*}};

\node[left] at (0,1)  {$ h_0 $};
\node[right] at (1,1)  {$ h_1 $};
\node[below=4pt, right] at (4,3)  {$ a = h_0 b_{20} $};
\node[left] at (8,6)  {$ a^2 $};
\node[below] at (8,4)  {$ b = b_{20}^2 $};
}

\draw[fill] (-10.3,7.7) circle [radius=0.1];
\node at (-8.5,7.7)  { \small{means $\M$} };
\node[mark size=1.5pt,color=red, right] at (-10.3-0.88,6.7) {\pgfuseplotmark{triangle*}} ;
\node at (-7.4,6.7)  { \small{means $\M/\tau \cong \F_2$} };
\node at (-7,5.2)  { \underline{(stem, filtration, weight):} };
\node at (-7,4)  {$|\tau|=(0,0,-1)$ };
\node at (-7,3)  {$h_0 = 2$ and $|h_0|=(0,1,0)$ };
\node at (-7,2)  {$h_1 = \eta$ and $|h_1|=(1,1,1)$ };
\node at (-7,1)  {$|a|=(4,3,2)$ };
\node at (-7,0)  {$|b|=(8,4,4)$ };
\end{tikzpicture}
\caption{The $E_{\infty}$-page of the Adams spectral sequence computing $\piss(\kq)$.} \label{fig:htpyofkO}
\end{figure}
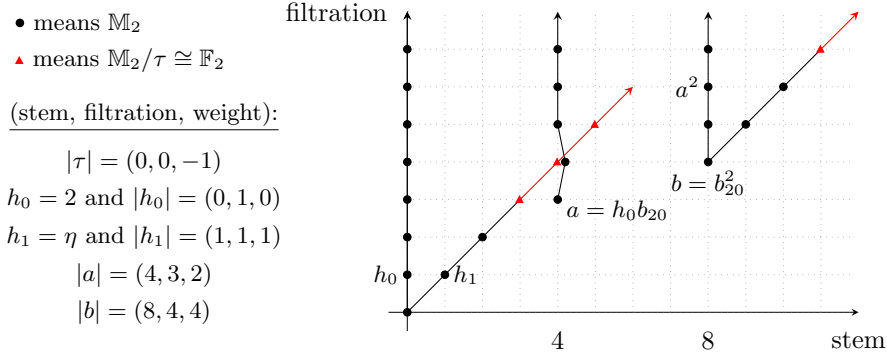

Consider now the hermitian $K$-theory spectrum $KQ$ defined in \cite{HornbostelHermitian} and studied in \cite{ROhermitian}. The paper \cite{IS} defines its connective cover $\kq$ over $\Spec \C$, by taking appropriate $C_2$-fixed points (although it is denoted by $ko$ in that paper). It also computes its coefficients and mod 2 homology
\begin{equation*}
\piss(\kq) \cong \quotient{\hat{\Z}_2[\tau, \eta, a,b]}{2 \eta, \tau \eta^3, a \eta, a^2 = 4b} \quad \text{ and } \quad \HFho( \kq) \cong \quotient{ \F_2[\tau][\xi_1^2,\xi_2, \ldots][\tau_2, \tau_3,\ldots]}{\tau_i^2 = \tau \xi_{i+1}}.
\end{equation*} 
To explain the homotopy ring $\piss(\kq)$, Figure \ref{fig:htpyofkO} displays the $E_{\infty}$-page of the motivic Adams spectral sequence computing $\piss(\kq)$. The horizontal axis represents the stem, i.e., the $s$ in $\pi_{s,w}(\kq)$, while the vertical axis represents the Adams filtration.
As it is usually done with motivic charts, the weight $w$ in $\pisw(\kq)$ is suppressed from the chart and one can imagine it on a third axis perpendicular to the page.

In this Section we consider the $\Ct$-induced spectra that we denote by
\begin{equation*}
\kglCt \coloneqq \kgl \smas \Ct \qquad \text{ and } \qquad \kqCt \coloneqq \kq \smas \Ct.
\end{equation*}
Both of them are $\Ct$-algebras, where $\kgl$ is an $E_{\infty}$ algebra as being the smash product of two $E_{\infty}$ rings.

\subsubsection*{The case of algebraic $K$-theory $\kglCt$} The fact that both its homotopy and homology are $\tau$-free makes the description of $\kglCt$ straightforward. Indeed, by Lemma \ref{lem:htpyhlgyoftaufree} we immediately get
\begin{equation*}
\piss( \kglCt) \cong \hat{\Z}_2[v_1] \qquad \text{ and } \qquad \HFho( \kglCt) \cong \F_2[\xi_1,\xi_2, \ldots] \otimes E(\tau_2, \tau_3,\ldots).
\end{equation*}

\vspace{0.5cm}

\subsubsection*{The case of hermitian $K$-theory $\kqCt$} Its homology is $\tau$-free and so again we immediately get
\begin{equation*}
\HFho(\kqCt) \cong \F_2[\xi_1^2,\xi_2, \ldots] \otimes E(\tau_2, \tau_3,\ldots).
\end{equation*}
Its homotopy is more interesting as it is not $\tau$-free, and we will get contributions both from the cokernel and kernel of multiplication by $\tau$. Moreover, a surprising fact occurs as there is a hidden extension which makes $\kqCt$ contain the periodicity element $v_1^2$ in its homotopy.

\begin{prop}
The homotopy ring $\piss(\kqCt)$ has the presentation
\begin{equation*}
\piss( \kqCt) \cong \quotient{ \hat{\Z}_2[\eta, v_1^2]}{2 \eta}.
\end{equation*}
\begin{proof}
The usual cofiber sequence \eqref{eq:cofseqtau} for $\Ct$, smashed with $\kq$ gives the cofiber sequence 
\begin{equation*}
\Sigma^{0,-1} \kq \stackrel{\tau}{\lto} \kq \stackrel{i}{\lto} \kqCt \stackrel{p}{\lto} \Sigma^{1,-1} \kq.
\end{equation*}
Since the homology $\HFho\left( \Sigma^{0,-1} \kq \right)$ is $\tau$-free, we ge  the short exact sequence
\begin{equation*}
0 \lto \HFho \left( \Sigma^{0,-1} \kq \right) \stackrel{\tau}{\lto} \HFho\left( \kq \right) \stackrel{i}{\lto} \HFho\left( \kqCt\right) \lto 0
\end{equation*}
in homology. For any motivic spectrum $X$, denote by $\Ext^{\ast}(X)$ the trigraded term
\begin{equation*}
\Ext_{\ACdual\text{-comod}}^{\ast,\ast,\ast}(\HFho(\Ss), \HFho(X))
\end{equation*}
that represents the $E_2$ page of the motivic Adams spectral sequence for $X$. We use the indicated grading in $\Ext^{\ast}(X)$ to denote the homological degree in $\Ext$, i.e., the Adams filtration on the $E_2$ page. From the above short exact sequence, we get a long exact sequence in $\Ext$-groups
\begin{equation*}
\cdots \stackrel{\tau}{\lto} \Ext^{\ast}(\kq) \stackrel{i_{\ast}}{\lto} \Ext^{\ast}(\kqCt) \stackrel{p_{\ast}}{\lto}  \Ext^{\ast+1}(\Sigma^{0,-1} \kq) \stackrel{\tau}{\lto} \cdots,
\end{equation*}
i.e., a long exact sequence in $E_2$ pages. This gives short exact sequences
\begin{equation*}
0 \lto \quotient{\Ext^{\ast}(\kq)}{\tau} \stackrel{i_{\ast}}{\lto} \Ext^{\ast}(\kqCt) \stackrel{p_{\ast}}{\lto}  {}_{\tau}\Ext^{\ast+1}(\Sigma^{0,-1} \kq) \lto 0,
\end{equation*}
where the left term is the cokernel of $\tau$ while the right term is the $\tau$-torsion. Since $i$ is a ring map, the term $\Ext(\kq)/\tau$ includes as a subring of $\Ext(\kqCt)$. However, this cokernel can act non-trivially on the $\tau$-torsion part, giving potential extension problems to solve. Since the motivic Adams spectral sequence for $\kq$ collapses at the $E_2$ page with no hidden extensions, the term $\Ext(\kq)$ is given by the Figure \ref{fig:htpyofkO} on page \pageref{fig:htpyofkO}. These two pieces assemble to give the additive description of the $E_2$ page of the motivic Adams spectral sequence for $\kqCt$ as described in Figure \ref{fig:htpykOCtadditively}.
\begin{figure}[h] 
\begin{tikzpicture}[scale=0.5] 
\centering

\draw [very thin, -stealth] (0,-0.5) grid (0,8);
\node [left= 4pt] at (0,8) { filtration };
\draw [very thin, -stealth] (-0.5,0) -- (12,0);
\node [below=4pt] at (12,0) { stem };
\draw[step=1cm,gray,very thin, dotted] (0,0) grid (11.9,7.9);
\node [below=4pt] at (4,0) { $4$ };
\node [below=4pt] at (8,0) { $8$ };

\small{
\foreach \x in {0,...,7} { \node[mark size=1.5pt,color=red, right] at  (0-0.88,\x) {\pgfuseplotmark{triangle*}};  } 
\draw [-stealth, red] (0,0) -- (0,8);
						   			   
\foreach \x in {5,...,7} { \node[mark size=1.5pt,color=red, right] at  (4-0.88,\x) {\pgfuseplotmark{triangle*}}; } 
\node[mark size=1.5pt,color=red, right] at (4-0.88,3) {\pgfuseplotmark{triangle*}};
\node[mark size=1.5pt,color=red, right] at (4.2-0.88,4) {\pgfuseplotmark{triangle*}};
\draw [red] (4,3) -- (4.2,4) -- (4,5);
\draw [-stealth, red] (4,5) -- (4,8);

\foreach \x in {4,...,7} { \node[mark size=1.5pt,color=red, right] at (8-0.88,\x) {\pgfuseplotmark{triangle*}}; } 
\draw [-stealth,red] (8,4) -- (8,8);

\node[mark size=1.5pt,color=red, right] at (1-0.88,1) {\pgfuseplotmark{triangle*}};
\node[mark size=1.5pt,color=red, right] at (2-0.88,2) {\pgfuseplotmark{triangle*}};
\draw [red] (0,0) -- (1,1)-- (2,2)-- (3,3);
\draw [-stealth, red] (3,3) -- (6,6);
\node[mark size=1.5pt,color=red, right] at (3-0.88,3) {\pgfuseplotmark{triangle*}};
\node[mark size=1.5pt,color=red, right] at (4-0.88,4) {\pgfuseplotmark{triangle*}};
\node[mark size=1.5pt,color=red, right] at (5-0.88,5) {\pgfuseplotmark{triangle*}};

\draw [red] (8,4) -- (11,7);
\draw [-stealth, red] (11,7) -- (12,8);
\node[mark size=1.5pt,color=red, right] at (9-0.88,5) {\pgfuseplotmark{triangle*}};
\node[mark size=1.5pt,color=red, right] at  (10-0.88,6) {\pgfuseplotmark{triangle*}};
\node[mark size=1.5pt,color=red, right] at (11-0.88,7) {\pgfuseplotmark{triangle*}};

\node[left] at (0,1)  {$ h_0 $};
\node[right] at (1,1)  {$ \overline{h}_1 $};
\node[below=4pt, right] at (4,3)  {$ \overline{a}$};
\node[left] at (8,6)  {$ \overline{a}^2 $};
\node[below] at (8,4)  {$ \overline{b} $};

\draw [-stealth, red] (4,2) -- (7,5);
\node[mark size=1.5pt,color=red, right] at (4-0.88,2) {\pgfuseplotmark{triangle*}};
\node[mark size=1.5pt,color=red, right] at (5-0.88,3) {\pgfuseplotmark{triangle*}};
\node[mark size=1.5pt,color=red, right] at (6-0.88,4) {\pgfuseplotmark{triangle*}};
\node[right, below] at (4,2)  {$ \widetilde{h_1^3} $};

\draw [-stealth, red] (12,6) -- (13,7);
\draw [-stealth, red] (12,7) -- (12,8);
\node[mark size=1.5pt,color=red, right] at (12-0.88,6) {\pgfuseplotmark{triangle*}};
\node[mark size=1.5pt,color=red, right] at (12-0.88,7) {\pgfuseplotmark{triangle*}};

\draw [densely dashed, thick, blue] (4,2) to [bend left=45] (4,3);
\draw [densely dashed, thick, blue] (12,6) to [bend left=45] (12,7);
}

\node[mark size=1.5pt,color=red, right] at (-10.7-0.88,7.7) {\pgfuseplotmark{triangle*}};
\node at (-7.4,7.7)  { \small{means $\M/\tau \cong \F_2$} };
\node at (-7,6.2)  { \underline{(stem, filration, weight):} };
\node at (-7,5)  {$h_0 = 2$ and $|h_0|=(0,1,0)$ };
\node at (-7,4)  {$\overline{h}_1 = \overline{\eta}$ and $|\overline{h}_1|=(1,1,1)$ };
\node at (-7,3)  {$|\overline{a}|=(4,3,2)$ };
\node at (-7,2)  {$|\overline{b}|=(8,4,4)$ };
\node at (-7,1)  {$\widetilde{h_1^3} = \widetilde{\eta^3}$ and $|\widetilde{h_1^3}|=(4,2,2)$ };
\end{tikzpicture}
\caption{The $E_2$ page of the motivic Adams spectral sequence for $\kqCt$ as an $\F_2$-vector space.} \label{fig:htpykOCtadditively}
\end{figure}
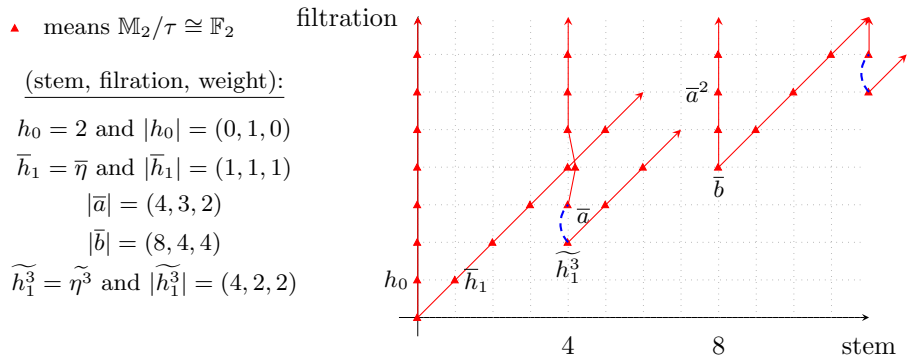
It still remains to solve the possible extension problems and possible Adams differentials. The only possible extension is whether or not $2 \cdot \widetilde{h_1^3} = \overline{a}$, as indicated in Figure \ref{fig:htpykOCtadditively}. Consider the Toda bracket $\langle \tau, \eta^3, 2 \rangle$ as in the diagram
\begin{center}
\begin{tikzpicture}
\matrix (m) [matrix of math nodes, row sep=3em, column sep=3em]
{ S^{3,2} & S^{3,2} & \Sigma^{0,-1} \kq & \kq \\
          &         & \Sigma^{-1,0} \kqCt & \\
          &         & \Sigma^{-1,0} \kq, & \\};
\path[thick, -stealth, font=\small]
(m-1-1) edge node[above] {$ 2 $} (m-1-2)
(m-1-2) edge node[above] {$ \eta^3 $}  (m-1-3)
(m-1-3) edge node[above] {$ \tau  $} (m-1-4)
(m-2-3) edge node[right] {$ p $} (m-1-3)
(m-3-3) edge node[right] {$ i $} (m-2-3);
\path[thick, -stealth, dashed, font=\small]
(m-1-2) edge node[below, left= 4pt] {$ \widetilde{\eta^3} $} (m-2-3)
(m-1-1) edge (m-3-3);
\end{tikzpicture}
\end{center}
where we have that $2 \cdot \widetilde{\eta^3} \in i_{\ast} \langle \tau, \eta^3, 2 \rangle$ by \cite[Section 3.1.1]{StableStems}. We can compute this bracket in the motivic May spectral sequence using May's convergence Theorem. See \cite{MayMasseyProd} for the original reference, and \cite[Theorem 2.2.3]{StableStems} for an exposition of the motivic version. More precisely, we can compute it on the motivic May $E_3$-page via the  differential $d_3(b_{20}) = \tau h_1^3$ (since $h_0h_1$ is already zero). This bracket has no indeterminacy giving
\begin{equation*}
\langle \tau, h_1^3, h_0 \rangle = \left\{ b_{20}h_0 \right\}.
\end{equation*}
Recall from Figure \ref{fig:htpyofkO} that $a = b_{20}h_0$ giving that indeed, in $\piss( \kqCt)$, there is an extension $2 \cdot \widetilde{h_1^3} = \overline{a}$. This $h_0$-extension appears as the round dotted line on Figure \ref{fig:htpykOCtadditively}. We now spell out the ring structure of this $E_2$ page. First observe that
\begin{equation*}
4 \left( \widetilde{h_1^3} \right)^2 = \left( 2 \widetilde{h_1^3} \right)^2 = \overline{a}^2 = 4 \overline{b}^2,
\end{equation*}
and because there are no possible extensions in that column, we get that $\left(\widetilde{h_1^3} \right)^2 = \overline{b}$. The $E_2$ page of the motivic Adams spectral sequence for $\kqCt$ has therefore the ring presentation 
\begin{equation*}
E_2 \cong \quotient{\F \left[ h_0, \overline{h}_1, \widetilde{h_1^3} \right] }{h_0 \overline{h}_1 }.
\end{equation*}
There are no possible Adams differentials on these 3 generators, and thus Figure \ref{fig:htpykOCtadditively} also represents the $E_{\infty}$ page of the Adams spectral sequence for $\kqCt$. Except the $h_0$-towers, there are no possible hidden extensions, giving the multiplicative description 
\begin{equation*}
\piss( \kqCt) \cong \quotient{\hat{\Z}_2[\eta, \widetilde{h_1^3}]}{2 \eta}.
\end{equation*}
Finally, we show that $\widetilde{h_1^3}$ detects the element $v_1^2$. We can smash the cofiber sequence
\begin{equation*}
\Sigma^{1,1} \kq \stackrel{\eta}{\lto} \kq \stackrel{i}{\lto} \kgl
\end{equation*}
with $\Ct$ to obtain the cofiber sequence
\begin{equation*}
\Sigma^{1,1} \kqCt \stackrel{\eta}{\lto} \kqCt \stackrel{\overline{i}}{\lto} \kglCt.
\end{equation*}
Since $i$ is a ring map, then so is the induced map $\overline{i}$. The ring map $\overline{i}$ sends the 8-fold Bott periodicity element $\overline{b} = \left( \widetilde{h_1^3} \right)^2$ to the 8-fold Bott periodicity element $v_1^4$, which forces $\widetilde{h_1^3}$ to be sent to $v_1^2$. The $E_2$ page of $\kqCt$ has therefore the ring presentation 
\begin{equation*}
\piss( \kqCt) \cong \quotient{\hat{\Z}_2[\eta, v_1^2]}{2 \eta}. \qedhere
\end{equation*}
\end{proof}
\end{prop}

\bibliographystyle{alpha}
\bibliography{mybibliography}

\end{document}